\documentclass[twocolumn]{autart}    

\usepackage{graphicx}          
\usepackage{xcolor}
\usepackage{amsmath}
\usepackage{amsfonts}
\usepackage{pgfplots}
\usepackage{tikz-network}
\pgfplotsset{compat=1.18} 
\usepackage{pifont}
\usepackage{url}
\usepackage{adjustbox}
\usepackage{booktabs}
\usepackage{multirow}
\usepackage{natbib}
\usepackage{mathtools}
\usepackage{algorithm}
\usepackage{subcaption}
\usepackage{hyperref}
\usepackage{stackengine}
\usetikzlibrary{arrows}
\usepackage[noend]{algpseudocode}
\usepackage{enumitem}
\usepackage{lipsum}

\usepgfplotslibrary{groupplots} 
\usepgfplotslibrary{statistics}
\usepgfplotslibrary{fillbetween}

\newcommand{\av}{\mathbf{a}}
\newcommand{\bv}{\mathbf{b}}
\newcommand{\cvec}{\mathbf{c}}
\newcommand{\ev}{\mathbf{e}}

\newcommand{\q}{\mathbf{q}}
\newcommand{\rv}{\mathbf{r}}

\newcommand{\uv}{\mathbf{u}}

\newcommand{\w}{\mathbf{w}}
\newcommand{\x}{\mathbf{x}}
\newcommand{\y}{\mathbf{y}}
\newcommand{\z}{\mathbf{z}}

\newcommand{\1}{\mathbf{1}}
\newcommand{\A}{\mathbf{A}}
\newcommand{\B}{\mathbf{B}}
\newcommand{\C}{\mathbf{C}}
\newcommand{\D}{\mathbf{D}}

\newcommand{\I}{\mathbf{I}}
\newcommand{\Lm}{\mathbf{L}}
\newcommand{\Pm}{\mathbf{P}}
\newcommand{\Q}{\mathbf{Q}}

\newcommand{\U}{\mathbf{U}}

\newcommand{\W}{\mathbf{W}}

\newcommand{\Y}{\mathbf{Y}}
\newcommand{\Z}{\mathbf{Z}}

\newcommand{\At}{\mathbfcal{A}}

\newcommand{\Ct}{\mathbfcal{C}}

\newcommand{\Qt}{\mathbfcal{Q}}
\newcommand{\Qts}{\mathbfcal{\scriptstyle{Q}}}

\newcommand{\Wt}{\mathbfcal{W}}

\newcommand{\Yt}{\mathbfcal{Y}}

\newcommand{\bigO}{\mathcal{O}}

\newcommand{\Rtensor}{\mathbb{R}^{I_1 \times I_2 \times \dots \times I_D}}

\newcommand{\Rcore}{\mathbb{R}^{R_{d} \times I_d \times R_{d+1}}}

\DeclareMathOperator{\vect}{vec}

\DeclareMathOperator{\VAF}{VAF}

\DeclareMathOperator*{\argmin}{arg\,min}
\DeclareMathAlphabet\mathbfcal{OMS}{cmsy}{b}{n}

\DeclarePairedDelimiterX{\inp}[2]{\langle}{\rangle}{#1, #2}

\RequirePackage{xcolor}
\definecolor{NiceRed}{RGB}{ 255, 92, 168}
\definecolor{NiceGreen}{RGB}{ 90, 168, 0}
\definecolor{NiceBlue}{RGB}{ 0, 152, 233}
\definecolor{NiceYellow}{RGB}{ 242, 147, 24}
\colorlet{NiceViolet}{NiceRed!50!NiceBlue}
\colorlet{NiceBrown}{NiceRed!50!NiceGreen}
\colorlet{NiceOrange}{NiceRed!50!NiceYellow}
\colorlet{NiceCyan}{NiceGreen!50!NiceBlue}

\definecolor{MatchingBlue}{RGB}{ 19, 123, 255}
\definecolor{MatchingPink}{RGB}{ 255, 6, 122}
\definecolor{MatchingCyan}{RGB}{19, 182, 148}
\definecolor{DetailOrange}{RGB}{235, 114, 70}

\definecolor{blue}{RGB}{68,119,170}
\definecolor{cyan}{RGB}{102,204,238}
\definecolor{green}{RGB}{34,136,51}
\definecolor{yellow}{RGB}{204,187,68}
\definecolor{red}{RGB}{238,102,119}
\definecolor{purple}{RGB}{170,51,119}
\definecolor{grey}{RGB}{187,187,187}
\definecolor{orange}{RGB}{238,119,51}

\begin{document}

\begin{frontmatter}

\title{Automatic Structure Identification for Highly Nonlinear MIMO Volterra Tensor Networks\thanksref{footnoteinfo}} 

\thanks[footnoteinfo]{Eva Memmel, and
thereby this work, is supported by the Delft University of
Technology AI Labs program.}
\author[Delft]{Eva Memmel}\ead{e.m.memmel@tudelft.nl},    
\author[Delft]{Kim Batselier}\ead{k.batselier@tudelft.nl},               

\address[Delft]{Delft University of Technology, Netherlands}  


\begin{abstract}                          
The Volterra Tensor Network lifts the curse of dimensionality for truncated, discrete times Volterra models, enabling scalable representation of highly nonlinear system.
This scalability comes at the cost of introducing randomness through initialization, and leaves open the challenge of how to efficiently determine the hyperparameters model order and memory length. 
In this paper, we present a unified framework that simultaneously addresses both challenges: We derive two algorithms that incrementally increase the model order and memory length along. Further we proof that the updates are performed along conjugate directions by establishing a mathematical equivalence between our proposed algorithms and equality constrained least squares systems. We present several strategies how to use our proposed algorithms for initialization and hyperparameter selection.
In numerical experiments, we demonstrate that our proposed algorithms are more accurate and efficient than the state-of-the-art Volterra Tensor Network and achieve competitive results to several state-of-the-art Volterra models.
\end{abstract}

\end{frontmatter}

\section{Introduction}
Nonlinear system identification is typically an iterative process that involves a wide range of user choices \cite{schoukens2019nonlinear}. 
In the context of black-box models, a challenge lies in balancing structural model errors against the model complexity \cite{schoukens2019nonlinear}. A popular black-box model class is the discrete time, truncated Volterra Series. 
For single-input single-output (SISO), it is given by {\small%
\vspace{-0.75cm}
\begin{align}
    \hat{y}(n) =\sum_{d = 0}^{D} \sum_{m_1,\dots,m_d}^{M-1} w_d(m_1,\dots ,m_d) \prod_{j=1}^d u(n-m_j),\label{eq:SISOVolterra}
\end{align}}%
where  $u(n),\ \hat{y}(n)$ is the input and output at the time step $n = 0,\dots,N-1$, respectively, and all parameters $w_d$ form the $d$th order Volterra kernel, $D$ denotes the order of the Volterra model and $M$ the memory length.
Equation (\ref{eq:SISOVolterra}) can be exptended to multiple-input, multiple-output (MIMO) with $P$ inputs and $L$ outputs.\\
A well-known limitation of the Volterra series is the Curse of Dimensionality, describing the exponential growth of the total number of model parameters $w_d$ with the model order $D$.
As a result, computational costs commonly restrict the feasible model structures to low complexity models. 
Unfortunately, those overly simple models often cause structural model errors due to the omission of higher order interaction terms (not high enough $D$) or the exclusion of particular past inputs (not high enough $M$) \cite{birpoutsoukis2018volterra, skyvulstad2023regularised, immordino2025parametric}.
Moreover, $D$ and $M$ are seldom known \textit{a priori} and therefore must be treated as hyperparameters, which are commonly tuned with cross-validation, grid search \cite{martinez2025tensor} or Bayesian optimization \cite{birpoutsoukis2018volterra}.
Unfortunately, these methods require training a new model order for each hyperparameter combination. 
Due to the Curse of Dimensionality, this is particularly computationally expensive for Volterra models. 
Typically, training one Volterra model incurs the computational costs $\bigO(M^{3D})$, increasing $D$ or $M$ costs $\bigO(M^{3(D+1)})$ or $\bigO((M+M_\Delta)^{3D})$, respectively. 
The costs of model structure selection can thus be reduced by reducing the trainings costs of a single model or the costs of the model order selection method.\\
\\
One method, that focuses on reducing training costs by lifting the Curse of Dimensionality is the low-rank Volterra Tensor Network (VTN) \cite{batselier2017tensor}.
The VTN uses the Tensor Train (TT) decomposition \cite{oseledets2011tensor} to compress all Volterra kernels at once into $D$ TT-cores of rank $R$. 
The train the model, the TT-cores are randomly initialized for a fixed $D$, $M$ and $R$, and then trained with the Alternating Linear Scheme (ALS) \cite{rohwedder2013local}.
The ALS updates one TT core at a time, moving back and forth across all TT-cores in $K$ sweeps. 
For a single model, the VTN reduces the storage complexity from $\bigO(M^D)$ to $\bigO(DMR^2)$ and the computational complexity to $\bigO(KD(M^3R^6))$\cite{batselier2017tensor}.
Although the VTN thereby enables the efficient training of highly nonlinear models, an extensive search of hyperparameters remains costly.
Due to the required random restarts, several VTN models need to be trained for each possible hyperparameter configuration. 
\\
\\
The main contribution of our paper is the derivation of an efficient method to select the hyperparameters $D$ and $M$ in VTN format.
To increase e.g. $D$ to $D+1$, we first reinterpret the least squares (LS) problem of a $D$th order Volterra model as an $D+1$ order equality constraint least squares (LSE) problem with zero constraints.
Then, we lift the quality constraints and solve the corresponding $D+1$th LS problem.
We prove, that this is equivalent to training a $D+1$th order residual model on the current residual $\rv_D$ of a $D$th order model.
Moreover, it ensures that the training outputs $\hat{\y}_D$ of the $D$th order model are orthogonal to the outputs $\hat{\y}_{\Delta} = \hat{\y}_{D+1} - \hat{\y}_{D}$ of the residual model.
We incorporate the LSE formulation in the VTN structure and provide guarantees for the training error after the increase step. 
Moreover, we demonstrate how to train the residual VTN with just a single core update, such that increasing $D$ to $D+1$ incurs the marginal costs of $\bigO(M^3R^6)$.
\\
\\
The resulting hybrid algorithm combines characteristics of forward regression \cite{billings1988identification, chen1989orthogonal}, gradient boosting \cite{friedman2001greedy} and the ALS steepest descend (ALS-SD) \cite{dolgov2014alternating}. 
The ALS-SD relies on a standard steepest descent step, which restricts the proposed method to increasing the rank $R$ for a given $D$ and $M$.
In contrast, we perform the steepest descent step in the prediction space, just as gradient boosting.
Instead of fixing previously trained parameters, we retrain them in the increase step to get closer to the symmetric minimum norm solution of the global LS problem.
Similar to the ALS-SD, we alternate between increase steps and standard ALS updates for improved accuracy.
This introduces the number of \textit{sweeps} performed between each increase step as a new hyperparameter.
While \textit{sweeps} is fixed in the ALS-SD, we demonstrate the trade-off between the model order $D$ and \textit{sweeps}. 
We show that, similar to the learning rate in gradient boosting, a lower number of \textit{sweeps} leads to a higher number of the optimal model order. 
\\
\\
We leverage the orthogonality of the outputs $\hat{\y}_D$ and $\hat{\y}_{\Delta}$ to evaluate the isolated contribution of the residual VTN on standard metrics.   Similar to the forward regression, these metrics informs the decision of whether to accept or reject the incremental increase. Whereas \cite{billings1988identification, chen1989orthogonal} impose a sparse model assumption, our proposed low-rank approach allows for more complex model structures.
Based on the evaluation of the incremental increase, we propose a method for automatic model structure selection. Starting with a low complexity model (small $D$ and $M$), at each increase step it is automatically decided with parameter to increase.
So, instead of performing a $D \times M$ grid search, this method finds a path through this grid with just one model.
As a welcome addition, our proposed method yields a deterministic initialization strategy for the VTN, making random restarts obsolete.
Instead of randomly initializing the TT-cores, we estimate a unique low-order (small $D$) VTN model as a convex least squares problem \cite{batselier2021enforcing} to increase $D$ until the desired model order is reached. 
\\
\\
In the experiments we showcase different approaches to select the hyperparameters $D$ and $M$.
We use the synthetical data set proposed in \cite{batselier2021enforcing} and the nonlinear benchmark Cascades \cite{schoukens2016cascaded} and compare to Volterra model in the literature. 
In our numerical experiments, we demonstrate that even simultaneously determining $D$ and $M$ achieves increased or competitive accuracy and a significant speed-up compared to the state-of-the-art (SOTA) VTN \cite{batselier2017tensor}.

\section{Preliminaries and Problem Setup}
We introduce our notation in Section \ref{sec:notation}, provide background on TNs in Section \ref{sec:basics_tensors} and present the problem setup based on the VTN \cite{batselier2017tensor} in Section \ref{sec:basics_VTN}.
We assume basic familiarity with TNs and use standard tensor notation defined in the literature. Please note, that tensor notation is more complex than matrix notation.
A comprehensive review on TNs is given e.g. in \cite{kolda2009tensor, oseledets2011tensor}, a tutorial on TNs for system identification in \cite{batselier2022low}. 
\subsection{Notation}\label{sec:notation}
A $d$\textsuperscript{th} dimensional array is a $d$\textsuperscript{th}-order tensor. 
We distinguish between scalars $a$ ($D=0$), vectors $\av$ ($D=1$), matrices $\A$ ($D=2$) and higher-order tensors $\At$ ($D \geq 3$). 
A scalar tensor element is accessed with a vertical bar combined with subscript indices in brackets, e.g. for a matrix $\A \in \mathbb{R}^{J \times K}$, $a|_{\left[j,k\right]}$ denotes its entry in the $j$th row and $k$th column; $\av|_{[j,:]}$ denotes its $j$th row. 
We denote a zero-vector (one-vector) in $\mathbb{R}^K$ as $\mathbf{0}_K$ ($\mathbf{1}_K$), respectively.
To extend those objects to higher orders, we specify the respective dimensions, e.g. as $\mathbf{0}_{J\times K}$.
We denote the $k$th standard basis vector of $\mathbb{R}^K$ as $\ev_{K,k}$, such that $\ev_{K,k}$ has the entry $1$ at the $k$-th position and $0$ elsewhere.
$\mathbf{I}_K$ denotes the identity matrix in $\mathbb{R}^{K\times K}$, $\Q$ indicates orthogonality. The Kronecker product is denoted by $\otimes$. 
The $d$-times repeated Kronecker product is defined as $\x^{\otimes D}:= \x \otimes \x \otimes \dots \otimes \x $.\\
\\
We distinguish different contexts for the parameter vector $\w$, which contains all Volterra Kernels up to order $D$ as follows:
We indicate the minimum-norm solution of any LS/LSE problem with $\hat{\w}$.
A solution obtained in VTN format is specified with an additional superscript as $\hat{\w}{}^{\text{TT}}$.
A single, vectorized TT-core is denoted by its core-index, e.g. $\w^{(d)}$.
The context of a standard LS problem is indicated by specifying the model order, e.g. $\w_D$;
the context of an LSE problem is indicated by specifying its particular formulation, e.g. $\w_\text{LSE}$.
\subsection{Tensors and Tensor Train Decomposition}\label{sec:basics_tensors}
Any tensor object can be reshaped into a vector or matrix and vice versa by rearranging its entries.
We define reshaping a $D$th order tensor $\At\in \Rtensor$ into a column vector $\av \in \mathbb{R}^{I_1\dots I_D}$ as $\av = \vect{(\At)}$.
\begin{defn}[Left- and Right-Unfolding \cite{rohwedder2013local}] \label{def:right_left_unfold}
    The left- and right-unfolding of a $3$rd order tensor $\At \in \Rcore$ are obtained by reshaping $\At$ into the matrix $\A|_{\text{left}} \in \mathbb{R}^{R_{d}I \times R_{d+1}}$ with rows indexed by $r_{d} + iR_{d}$ and $\A|_{\text{right}} \in \mathbb{R}^{R_{d} \times IR_{d+1}}$ with columns indexed by $i + r_{d+1}I$, respectively.
\end{defn}
Similar to matrices, tensor are multiplied by summing over a shared index. 
This is also referred to as contraction or $k$-mode product, defined in Appendix \ref{app:properties_and_definitions}. 
Tensor networks (TNs), or tensor decompositions, generalize matrix decompositions to higher orders.
This paper uses the tensor train (TT) decomposition \cite{oseledets2011tensor}. 
The multilinearity of the TT decomposition yields the definition of the frame matrix \cite{dolgov2014alternating}. 
\begin{defn}[Tensor Train \cite{oseledets2011tensor}] \label{def:TT}
The rank \\($R_1,\dots, R_{D+1}$) TT-decomposition with $R_1 = R_{D+1}=1$ decomposes a $D$th order tensor $\At\in \Rtensor$ in a set of $D$ $3$rd order tensors $\At^{(1)},\dots,\At^{(D)}$, such that 
\begin{align}
    a\big|_{\left[i_1,i_2,\dots,i_D\right]} &=\sum_{r_1=1}^{R_1}\dots\sum_{r_{D+1}=1}^{R_{D+1}} \prod_{d=1}^D a^{(d)}\big|_{[r_d,\,i_d,\,r_{d+1}]}.\label{eq:TTdefinition}
\end{align}
All tensors $\At^{(d)} \in \Rcore$ are called TT-cores. The operator $\mathbf{\mathrm{TT}_{D}}(\cdot,\dots,\cdot):$
$\mathbb{R}^{R_1 \times I_1\times R_2}\times \dots\times \mathbb{R}^{R_D\times I_D \times R_{D+1}}\rightarrow \mathbb{R}^{I_1\dots I_D}$ \begin{align}
    \av = \vect(\At) = \mathbf{\mathrm{TT}_{D}}(\At^{(1)},\dots,\At^{(2)})
\end{align}
undoes the TT decomposition of $\At$ by contracting over all TT-cores and vectorizes the resulting tensor.
\end{defn}
\begin{defn}[Frame Matrix \cite{dolgov2014alternating}]
TT is linear in each TT-core, such that $\vect \left(\At \right) = \A_{\setminus d} \vect \left( \At^{(d)}\right)$, where $\A_{\setminus d} \in \mathbb{R}^{I_1\dots I_D \times R_d I_d R_{d+1}}$ is called frame matrix . 
The frame matrix is computed as 
\begin{align}
    \A_{\setminus d} = \A_{k>d} \otimes \I_{I_d} \otimes \A_{k<d}^\top, \label{eq:framematrix}
\end{align}
where $\A_{k<d} \in \mathbb{R}^{\I_1\dots I_{d-1} \times R_{d-1}},\ (\A_{k>d} \in \mathbb{R}^{R_d \times I_{d+1}\dots I_D})$ denotes the matrices resulting of the contraction of all TT-cores left (right) of the $d$\textsuperscript{th} core and $\otimes$ is the Kronecker product. 
\end{defn}
\begin{defn}[Left- and right-orthogonal \cite{rohwedder2013local}]
    A matrix $\Q \in \mathbb{R}^{J,K}$ is left-, or right-orthogonal, if $\Q^\top \Q = \I_K$, resp. $\Q \Q^\top = \I_J$.
    Accordingly, a TT-core $\Qt^{(d)} \in \Rcore$ is left-, or right-orthogonal, if
    \begin{align}
        \Q\big|_{\text{left}}^\top\Q\big|_{\text{left}} = \I_{R_{d+1}}, \text{ resp. } \Q\big|_{\text{right}}\Q\big|_{\text{right}}^\top = \I_{R_{d}}.
    \end{align}
\end{defn}
\begin{defn}[Site-$d$-mixed-canonical form \cite{dolgov2014alternating}] 
    A TT is in site-$d$-mixed-canonical form, if the frame matrix $\Q_{\setminus d}$ is left-orthogonal.
\end{defn} 
The frame matrix $\Q_{\setminus d}$ is left-orthogonal, if the matrices $\Q_{p>d}$ and $\Q_{p<d}^\top$ are left- or right-orthogonal, i.e. the corresponding TT-cores are left- or right orthogonal. \cite{dolgov2014alternating}.
To orthogonalize $\Q_{\setminus d}$\cite{dolgov2014alternating}, we start from the first ($D$\textsuperscript{th}) core. A QR- (LQ-)decomposition is computed of the left- (right-)unfolding of the $k$th TT-core.
The orthogonal part is reshaped back into a TT-core, the upper- (lower-) triangular part is moved to the $k+1$\textsuperscript{th} ($k-1$\textsuperscript{th}) TT-core. 
This process is repeated until the TT is in site-$d$-mixed-canonical form. 
The side-$d$-mixed-canonical form is used in the ALS \cite{rohwedder2013local} to ensure numerical stability.
\subsection{Problem Setup: Volterra Tensor Network}\label{sec:basics_VTN}
We introduce the SOTA VTN \cite{batselier2017tensor}, focusing on SISO system.
The extension from SISO to MIMO does not affect our proposed method and is briefly explained in \ref{rem:SISO2MIMO}.
The goal is to first reformulate Equation (\ref{eq:SISOVolterra}) as an ordinary least squares problem. 
Then, the inputs and Volterra Kernels are expressed as a TT, allowing to efficiently train the VTN by updating one TT-core at a time.
The SOTA VTN constructs a suitable matrix-expression for the inputs, followed by its TT-representation. 
The vector $\uv_n$ contains all $M$ lagged inputs and is defined as  
\begin{align}
\uv_n := \left(1 \ u(n) \dots u(n-M+1)\right)^\top \in \mathbb{R}^{I}, \label{eq:X_n}
\end{align}
where $I := M+1$. Next, $\uv_n^{\otimes_{D}}$ is defined as the $D$-times Kronecker product of $\uv_n$, such that
\begin{align}
 \uv_n^{\otimes_{D}} := \overbrace{\uv_n\otimes \uv_n \otimes \dots \otimes \uv_n}^{D\text{ times }} \in \mathbb{R}^{I^D}.\label{eq:KronProdStructure}
\end{align}
The matrix $\U^{\otimes_{D}} \in N\times I^D$ is formed by stacking Equation (\ref{eq:KronProdStructure}) for all $N$, such that
\begin{align}
    \U^{\otimes_{D}} = ( \uv_0^{\otimes_{D}}  \uv_1^{\otimes_{D}} \dots  \uv_{N-1}^{\otimes_{D}})^\top. \label{eq:input_TT}
\end{align}
Specifically, the Kronecker-product structure allows $\U^{\otimes_{D}}$ to be treated as a rank-$N$ TT, such that we refer to Equation (\ref{eq:input_TT}) also as \textbf{input-TT}. 
In contrast to standard TTs, where ranks $R_1 = R_{D+1} = 1$, we have $R_1 = N$.
With the input-TT, the predicted outputs $\hat{\y}_D \in \mathbb{R}^{N}$ are expressed as 
\begin{align}
    \hat{\y}_D := \U^{\otimes_{D}}\hat{\w}_D,\label{eq:MIMOVolterraAllN}
\end{align} 
where $\hat{\w}_{D} \in \mathbb{R}^{I^D}$ contains the LS solution of all Volterra Kernels up to order $D$.
Directly solving for $\hat{\w}_{D}$ in Equation (\ref{eq:MIMOVolterraAllN}) requires $\bigO(I^{3D})$ flops, which can quickly become computationally infeasible for large $D$. 
The SOTA VTN model therefore approximates $\w_{D}$ by imposing a low-rank TT structure on $\w_{D}$.
In site $d$-mixed canonical form, the weight-TT $\w_D^{\text{TT}}$ is then given as
\begin{align}
\w_D^{\text{TT}} = \mathrm{TT}_{D}\left( \scriptstyle \Qt^{(1)},\dots,\Qt^{(d-1)},\Wt^{(d)}_D,\Qt^{(d+1)},\dots,\Qt^{(D)} \right). \label{eq:TT_level_0}
\end{align}
The SOTA VTN now solves the linear system in Equation (\ref{eq:MIMOVolterraAllN}) by updating one TT-core at a time, using the ALS.
In the \textbf{ALS-optimization step}, the $d$\textsuperscript{th} core of the weight-TT is updated by solving the LS problem
\begin{align}
\hat{\w}_D^{(d)}\, = \, \argmin_{\w_{D}^{(d)}} \, || \y\, - \,\U^{\otimes_{D}}\, \hat{\Q}_{\setminus d} \, \w_{D}^{(d)} ||,
 \label{eq:update_oneCore}
\end{align}
with the frame matrix $\hat{\Q}_{\setminus d}$ and $\hat{\w}_{D}^{(d)} = \vect (\widehat{\Wt}{}_{D}^{(d)})$.
Assuming the sequence of updated TT-cores goes from $d$ to $d-1$, the \textbf{ALS-orthogonalization step} brings the weight-TT into site-$d - 1$-mixed canonical form.
As the remaining cores are already orthogonalized, it is sufficient to perform the above-discussed orthogonalization process only for the $d$\textsuperscript{th} core.
The vector $\hat{\w}{}_{D}^{\text{TT}}$ and matrix $\U^{\otimes D}$ are never constructed explicitly. 
Instead, all computations are performed efficiently in TT-format.
\begin{rem}\label{rem:SISO2MIMO}
    To extend from SISO to MIMO, the $P$ input channels are stacked and transform $u(n)$ in Equation (\ref{eq:X_n}) into a vector $\uv(n)$ of length $P$, such that $I := PM+1$. The outputs and parameters in Equation (\ref{eq:MIMOVolterraAllN}) are stacked into the matrices $\hat{\Y}_{D} \ \in \ \mathbb{R}^{N \times L}$ and $\W_{D} \ \in \ \mathbb{R}^{I^D \times L}$, respectively. Consequently, $R_1=L$ in Equation (\ref{eq:TT_level_0}).
\end{rem}

\section{Incremental Increase of Hyperparameters}
We focus the derivation on the increase from $D$ to $D+1$, since it is mostly equivalent to the increase from $M$ to $M + M_{\Delta}$.
We briefly discuss the differences for the incremental increase of $M$ at the end of each subsection.
The extension to any $D$ (and $M$) follows by induction. 
\textbf{All proofs are given in the Appendix.}
\\
\\
Figure \ref{fig:outline_proof_LSE} shows the outline of the theoretical considerations in Section \ref{sec:LSE_VM}.
We solve two LS problems: one for a full $D$th-order Volterra model given in Equation (\ref{eq:SISOVolterra}) and one for the corresponding $D+1$th-order Volterra model. 
By reformulation the $D$th-order LS problem as a LSE problem, we derive the oblique projection operator $\Pm$, which maps the solution of the $D+1$th-order LS problem to the solution to the $D$th-order LS problem. 
The oblique projection operator then enables us to demonstrate two properties:
First, that the parameter $\hat{\w}_D$ are updated to $\hat{\w}_{D+1}$ in a conjugate direction, denoted as $\hat{\w}_{\Delta}$,
This conjugacy enables us to proof the second property, i.e. the orthogonality $\hat{\y}_{\Delta} \perp \hat{\y}_{D}$ and that $\hat{\w}_{\Delta}$ models the residual $\rv_D := \y - \hat{\y}_D$.
In Section \ref{sec:LSE_VTN}, we transfer the LSE formulation to the VTN format to perfrom all computations efficiently, while ensuring that properties shown in Section \ref{sec:LSE_VM} are preserved. 
The corresponding steps are depicted in Figure \ref{fig:incD_TensorDiagram} using a graphical notation, where each circle denotes a TT-core, each edge an index of a TT-core and connected edges indicate summations over the involved indices. 
\subsection{Increasing $D$ and $M$ via LSE Formulation }\label{sec:LSE_VM}
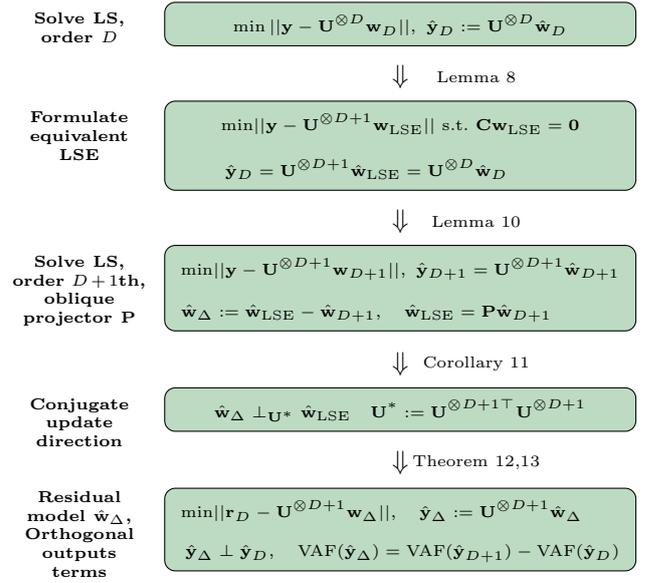
\begin{figure}
     \vspace{0.2cm}
    \centering
    \begin{minipage}{.45\textwidth}
        \hspace{-0.02\linewidth}%
        \begin{minipage}[t]{0.22\linewidth}
            \centering
            \tiny
            \textbf{Solve LS, order $D$}
        \end{minipage}%
        \hspace{0.03\linewidth}%
        \begin{minipage}[t]{0.75\linewidth}
            \centering \vspace{-0.75\baselineskip}
            \centering
            \begin{tikzpicture}
                \node[rectangle, rounded corners, draw=black, fill=green, fill opacity=0.3, text opacity = 1,minimum width=6.2cm](box) {
                \begin{minipage}{0.5\textwidth}
                    \tiny
                    \vspace{-0.4cm}
                        \begin{align*}
                            \min ||\y-\U^{\otimes D}\w_{D}||,\
                            \hat{\y}_{D} := \U^{\otimes D}\hat{\w}_{D}
                        \end{align*}
                    \end{minipage}
                    };
                    \node[below of=box, yshift=0.3cm] (arrow) {$\Downarrow$};
                    \node[right of=arrow] {\tiny Lemma \ref{lem:LSU_VM}};
            \end{tikzpicture}
        \end{minipage}%
    \end{minipage}
    \hfill
    \centering
    \begin{minipage}{0.45\textwidth}
        \vspace{0.15cm}
        \hspace{-0.02\linewidth}%
        \begin{minipage}[t]{0.22\linewidth}
            \centering
            \tiny
            \textbf{Formulate equivalent LSE}
        \end{minipage}%
        \hspace{0.03\linewidth}%
        \begin{minipage}[t]{0.75\linewidth}
            \centering \vspace{-0.75\baselineskip}
            \centering
            \begin{tikzpicture}
                \node[rectangle, rounded corners, draw=black, fill=green, fill opacity=0.3, text opacity = 1, minimum width=6.2cm] (box) {
                    \begin{minipage}{0.5\textwidth}
                    \tiny
                    \vspace{-0.4cm}
                        \begin{align*}
                            \min &||\y-\U^{\otimes D+1}\w_{\text{LSE}}|| \ \text{s.t.} \ \C \w_{\text{LSE}} = \mathbf{0} \\[0.2cm]
                            \hat{\y}_{D} &= \U^{\otimes D+1}\hat{\w}_{\text{LSE}} = \U^{\otimes D}\hat{\w}_{D}
                        \end{align*}
                    \end{minipage}
                    }; 
                \node[below of=box] (arrow) {$\Downarrow$};
                \node[right of=arrow] {\tiny Lemma \ref{lem:LSU_LSD1_VM}};
            \end{tikzpicture}
        \end{minipage}%
    \end{minipage}%
    \hfill%
    \begin{minipage}{0.45\textwidth}
        \vspace{0.15cm}
        \hspace{-0.02\linewidth}%
         \begin{minipage}[t]{0.22\linewidth}
            \centering
            \tiny
            \textbf{Solve LS, order $D+1$th, oblique projector $\Pm$}
        \end{minipage}%
        \hspace{0.03\linewidth}%
        \begin{minipage}[t]{0.75\linewidth}
            \centering \vspace{-0.75\baselineskip}
            \begin{tikzpicture}
                \node[rectangle, rounded corners, draw=black, fill=green, fill opacity=0.3, text opacity = 1,minimum width=6.2cm] (box) {
                    \begin{minipage}{0.5\textwidth}
                    \tiny
                    \vspace{-0.4cm}
                        \begin{align*}
                            \min &||\y-\U^{\otimes D+1}\w_{D+1}||,\ \hat{\y}_{D+1} = \U^{\otimes D+1}\hat{\w}_{D+1}\\[0.2cm]
                            \hat{\w}_{\Delta} &:= \hat{\w}_{\text{LSE}} - \hat{\w}_{D+1}, \quad \hat{\w}_{\text{LSE}} = \Pm\hat{\w}_{D+1} 
                        \end{align*}
                    \end{minipage}
                    }; 
                \node[below of=box, yshift=0.0cm] (arrow) {$\Downarrow$};
                \node[right of=arrow] {\tiny Corollary \ref{cor:conjugateUpdate_VM}};
            \end{tikzpicture}
        \end{minipage}%
    \end{minipage}%
\hfill%
    \begin{minipage}{0.45\textwidth}
        \vspace{0.15cm}
        \hspace{-0.02\linewidth}%
         \begin{minipage}[t]{0.22\linewidth}
            \centering
            \tiny
            \textbf{Conjugate update direction}
        \end{minipage}%
        \hspace{0.03\linewidth}%
        \begin{minipage}[t]{0.75\linewidth}
            \centering \vspace{-0.75\baselineskip}
            \begin{tikzpicture}
                \node[rectangle, rounded corners, draw=black, fill=green, fill opacity=0.3, text opacity = 1,minimum width=6.2cm] (box){
                    \begin{minipage}{0.5\textwidth}
                    \tiny
                    \vspace{-0.4cm}
                        \begin{align*}
                       \hat{\w}_{\Delta} \perp_{\U^*}\hat{\w}_{\text{LSE}} \quad \U^* := \U^{\otimes D+1 \top}\U^{\otimes D+1}            
                        \end{align*}  
                    \end{minipage}
                    }; 
                    \node[below of=box, yshift=0.3cm] (arrow) {$\Downarrow$};
                    \node[right of=arrow] {\tiny Theorem \ref{thm:residualModel_VM},\ref{thm:orthogonalUpdate_VM}};         
            \end{tikzpicture}
        \end{minipage}%
    \end{minipage}%
    \hfill%
    \begin{minipage}{0.45\textwidth}
        \vspace{0.075cm}
        \hspace{-0.02\linewidth}%
         \begin{minipage}[t]{0.22\linewidth}
            \centering
            \tiny
            \textbf{Residual model $\hat{\w}_\Delta$, Orthogonal outputs terms}
        \end{minipage}%
        \hspace{0.03\linewidth}%
        \begin{minipage}[t]{0.75\linewidth}
            \centering \vspace{-0.5\baselineskip}
            \begin{tikzpicture}
                \node[rectangle, rounded corners, draw=black, fill=green, fill opacity=0.3, text opacity = 1,minimum width=6.2cm] (box){
                    \begin{minipage}{0.5\textwidth}
                    \tiny
                    \vspace{-0.4cm}
                        \begin{align*}
                             \min &||\rv_D-\U^{\otimes D+1}\w_{\Delta}||, \quad \hat{\y}_{\Delta} := \U^{\otimes D+1}\hat{\w}_{\Delta}  \\[0.2cm]
                             \hat{\y}_{\Delta} &\perp \hat{\y}_{D}, \quad \VAF(\hat{\y}_{\Delta}) = \VAF(\hat{\y}_{D+1}) - \VAF(\hat{\y}_D)
                        \end{align*}
                    \end{minipage}
                    }; 
            \end{tikzpicture}
        \end{minipage}%
    \end{minipage}%
    \caption{Increasing the order from $D$ to $D+1$ via LSE. }\label{fig:outline_proof_LSE}
\end{figure}
We begin by providing relevant equations and definitions. The $D$th-order LS problem is given by
\begin{align}
   \min_{\w_{D}}||\y - \U^{\otimes D}\w_{D}||,\label{eq:VM_LS_D}
\end{align}
with $\w_D \in \mathbb{R}^{I^D}$. Analogous, the LS problem of a $D+1$th order Volterra model is given by
\begin{align}
   \min_{\w_{D+1}} ||\y - \U^{\otimes D+1}\w_{D+1}||, \label{eq:VM_LS_D1}
\end{align}
with $\w_{D+1} \in \mathbb{R}^{I^{D+1}}$.
Next, we formulate a $D+1$th-order LSE problem equivalent to the $D$th-order LS problem in Equation (\ref{eq:VM_LS_D}).
We achieve this equivalence by defining constraints $\C$ that zero out all additional parameters, which are associated with the $I^{D+1}-I^D = MI^D$ additional input terms created by expanding $\U^{\otimes D}$ to $\U^{\otimes D+1}$.
This $D+1$th-order LSE problem is given by 
\begin{align}
      \min_{\w_{\text{LSE}}} ||\y - \U^{\otimes D+1}\w_{\text{LSE}}|| \ \text{s.t.} \ \C  \w_{\text{LSE}} = \mathbf{0}_{MI^D}, \label{eq:VM_LSE}
\end{align}
where $\w_\text{LSE} \in \mathbb{R}^{I^{D+1}}$ represents a constrained $D+1$th-order Volterra model and $\C \in \mathbb{R}^{MI^D\times I^{D+1}}$ denotes the constraint matrix. 
We achieve the desired constraints by composing the rows $\C$ of $MI^D$ standard basis vectors $\ev_{I^{D+1},m}^\top$. We explicitly derive $\C$ in Appendix \ref{app:constrain_matrix_VM}.
Finally, we define the vectors
\begin{align}
    \hat{\w}_{\Delta} &:= \hat{\w}_{D+1} - \hat{\w}_{\text{LSE}}, \quad &\hat{\w}_{\Delta} &\in \mathbb{R}^{I^{D+1}}, \label{eq:updateVector_VM} \\
    \hat{\y}_{\Delta} &:= \U^{\otimes D+1} \hat{\w}_{\Delta}, \quad &\hat{\y}_{D+1} &\in \mathbb{R}^N, \label{eq:updateVector_outputs_VM}
\end{align}
such that $\hat{\w}_{\Delta}$ contains the updates from $\hat{\w}_{\text{LSE}}$ to $\hat{\w}_{D+1}$, and $\hat{\y}_{\Delta}$ denotes the outputs associated with $\hat{\w}_{\Delta}$.
\begin{rem}
    Lifting the constraints $\C$ in Equation (\ref{eq:VM_LSE}) yields the $D+1$th order LS problem in Equation (\ref{eq:VM_LS_D1}).
\end{rem}
We eliminate the constraints in Equation (\ref{eq:VM_LSE}) via the nullspace method \cite{bjorck1996numerical} to obtain an equivalent unconstrained LS (LSU) problem of lower dimension.
The LSU problem connects the $D$th-order LS problem in Equation (\ref{eq:VM_LS_D}) to the $D+1$th-order LSE problem in Equation (\ref{eq:VM_LSE}).
We then prove that $\hat{\w}_{D+1}$ increases the model order along conjugate directions to obtain the desired orthogonality $\hat{\y}_{\Delta} \perp \hat{\y}_{D}$. 
We show that $\hat{\w}_\Delta$ represents a model trained on the residual $\rv_D$.\\
\\
We obtain the LSU problem via a matrix by constructing a matrix $\Z \in \mathbb{R}^{I^{D+1}\times I^D}$, whose columns span an orthonormal basis of the nullspace $\mathcal{N}(\C)$ of $\C$. 
We derive in Appendix \ref{app:nullspace_matrix_VM}, that the columns of $\Z$ as $I^D$ consist of the standard basis vectors $\ev_{I^{D+1},i}$.
By re-parametrizing $\w_{\text{LSE}} = \Z \w_{\text{LSU}} $
we formulate the LSU problem as 
\begin{align}
   \min_{\w_{\text{LSU}}} ||\y - \U^{\otimes D+1}\Z\w_{\text{LSU}}||, \label{eq:VM_LSU}
\end{align}
with $\w_{\text{LSU}} \in \mathbb{R}^{I^D}$.
In Equation (\ref{eq:VM_LSU}), we can apply $\Z$ either to $\w_{\text{LSU}}$, or to $\U^{\otimes D+1}$. 
We summarize the respective effects in Lemma  \ref{lem:LSU_VM} and Corollary \ref{cor:identical_outputs}.
\begin{lem}\label{lem:LSU_VM}
    The solution $\hat{\w}_{\text{LSE}}$ can be formed by appending zero's to $\hat{\w}_{\text{LSU}}$. Moreover, $\Z^\top\hat{\w}_{\text{LSE}} = \hat{\w}_{\text{LSU}} = \hat{\w}_{D}$.
\end{lem}
\begin{cor}\label{cor:identical_outputs}
    The solutions $\hat{\w}_{\text{LSU}}$, $\hat{\w}_{\text{LSE}}$ and $\hat{\w}_{D}$ induce identical input-output behavior $\hat{\y}_D = \U^{\otimes D} \hat{\w}_{D}$.
\end{cor}
The constraints of the LSE problem in Equation (\ref{eq:VM_LSE}) effectively restrict the solution space of the $D+1$th-order LS problem given in Equation (\ref{eq:VM_LS_D}) in such a way, that each solution simultaneously fulfills the LS objectives in Equation (\ref{eq:VM_LS_D}) and the constraints in Equation (\ref{eq:VM_LSE}).
This relationship allows us to express $\hat{\w}_{\text{LSE}}$ in terms of $\hat{\w}_{D+1}$.
\begin{lem}\label{lem:LSU_LSD1_VM}
    The solution $\hat{\w}_{\text{LSE}}$ is an oblique projection of $\hat{\w}_{D+1}$, defined by the projector $\Pm$.
   With $(\cdot)^\dagger$ denoting the Moore-Penrose pseudoinverse, $\hat{\w}_{\text{LSE}}$ is given by
    \begin{align*}
        \hat{\w}_{\text{LSE}} = \underbrace{\Z\left(\U^{\otimes D+1}\Z\right)^\dagger \U^{\otimes D+1}}_{=:\ \Pm}\hat{\w}_{D+1}.
    \end{align*}
\end{lem}
In fact, any oblique projector $\Pm$ acts as an orthogonal projector with respect to the weighted inner product defined as $\langle \av, \bv\rangle_{\U^*} := \av^\top \U^* \bv$ \cite{stewart2011numerical}.
For the projector $\Pm$ defined in Lemma  \ref{lem:LSU_LSD1_VM}, $\U^* := \U^{\otimes D+1\top}\U^{\otimes D+1}$. 
This immediately yields the following Corollaries. 
\begin{cor}\label{cor:conjugateUpdate_VM}
    The vector $\hat{\w}_{\Delta}$ increases the order $D$ to $D+1$ along conjugate (or $\U^{\otimes D+1\top}\U^{\otimes D+1}$-orthogonal) update directions. Consequently,\begin{align*}
    \hat{\w}_{D+1}^\top \U^{\otimes D+1\top}\U^{\otimes D+1} \hat{\w}_{\text{LSE}} = 0.
    \end{align*} 
\end{cor}
Combining the results established above, we complete the derivation with the following theorems.
\begin{thm}\label{thm:residualModel_VM}
The parameters $\hat{\w}_{\Delta}$ model the residual $\rv_D = \y - \hat{\y}_D$, such that $\hat{\w}_{\Delta}$ is the solution to the minimization problem
\begin{align}
   \min_{\w_{\Delta}} ||\rv_D - \U^{\otimes D+1}\w_{\Delta}||. \label{eq:redisual_VM}
\end{align}
\end{thm}
\begin{thm}\label{thm:orthogonalUpdate_VM}
    The predicted $D$th-order outputs $\hat{\y}_{D}$ are orthogonal to the predicted outputs $\hat{\y}_{\Delta}$ in the standard Euclidian inner product space.
    We denote this orthogonality as $\hat{\y}_{\Delta} \perp \hat{\y}_{D}$. Moreover, $\hat{\y}_{D+1} = \hat{\y}_{D} + \hat{\y}_{\Delta}$.
\end{thm}
The derivation to increase $M$ is largely analogue. The main difference is that we do not expand the matrix $\U^{\otimes D}$ by adding another Kronecker product with $\uv_n$ to $\uv_n^{\otimes D}$.
Instead, we expand the matrix $\U^{\otimes D}$ by expanding each factor $\uv_n$ in Equation with (\ref{eq:KronProdStructure}) $M_\Delta$ additional lagged input terms.
The construction of $\C_M$ and $\Z_M$ are detailed in Appendix \ref{app:constrain_matrix_VM} and Appendix \ref{app:nullspace_matrix_VM}, respectively.
\subsection{Increasing $D$ and $M$ Efficiently in VTN-format}\label{sec:LSE_VTN}
\begin{figure}
    \centering
\SetEdgeStyle[TextFillOpacity = 0]
\begin{tikzpicture}

\begin{scope}[scale=0.85]
\Vertex[x=0,y=0,opacity=0, size=0.5]{1} 
\Vertex[x=0,y=-0.6,style={draw=none, fill=none}, size=.1]{W1}
\Edge[lw=1pt, position={left=1mm}](1)(W1)

\Vertex[x=0,y=-1.1,style={draw=none, fill=none}, size=.1]{W2}
\Edge[style={dotted}](W1)(W2)

\Vertex[x=0,y=-1.7,opacity=0, size=0.5]{3}
\Edge[lw=1pt, position={left=1mm}](W2)(3)

\Vertex[x=0,y=-3.5,style={fill=blue,fill opacity=0.3},position= below right, distance=-0.2cm,label=\tiny{\textcolor{blue}{$\widehat{\Wt}{}^{(d)}_D$}}, size=0.5]{5}
\Edge[lw=1pt, position={left=1mm}](3)(5)

\Vertex[x=0,y=-4.1,style={draw=none, fill=none}, size=.1]{W6}
\Edge[lw=1pt,position={left=1mm}](5)(W6)

\Vertex[x=0,y=-4.6,style={draw=none, fill=none}, size=.1]{W7}
\Edge[style={dotted}](W6)(W7)

\Vertex[x=0,y=-5.2,opacity=0,size=0.5]{8}
\Edge[lw=1pt, position={left=1mm}](W7)(8)

\Vertex[x=-0.75,y=-0,style={draw=none, fill=none}, size=.1]{I1}
\Vertex[x=-0.75,y=-1.7,style={draw=none, fill=none}, size=.1]{I3}
\Vertex[x=-0.75,y=-3.5,style={draw=none, fill=none}, size=.1]{I5}
\Vertex[x=-0.75,y=-5.2,style={draw=none, fill=none}, size=.1]{I8}

\Edge[lw=1pt, position={left=1mm}](I1)(1)
\Edge[lw=1pt, position={left=1mm}](I3)(3)
\Edge[lw=1pt, position={left=1mm}](I5)(5)
\Edge[lw=1pt, position={left=1mm}](I8)(8)

\node[rectangle, rounded corners, draw=none, fill=none, fill opacity=0.3, text opacity = 1, below of=8, align=center] {\rotatebox{90}{$\,=$} \\ \tiny $\hat{\w}_D^{\text{TT}}$};
\end{scope}

\begin{scope}[xshift=1.6cm, scale = 0.85]
\Vertex[x=0,y=0,opacity=0, size=0.5]{1} 
\Vertex[x=0,y=-0.6,style={draw=none, fill=none}, size=.01]{W1}
\Edge[lw=1pt, position={left=1mm}](1)(W1)

\Vertex[x=0,y=-1.1,style={draw=none,fill=none}, size=.1]{W2}
\Edge[style={dotted}](W1)(W2)

\Vertex[x=0,y=-1.7,opacity=0, size=0.5]{3}
\Edge[lw=1pt, position={left=1mm}](W2)(3)

\Vertex[x=0,y=-2.6,opacity=0, size=0.5, position= above right, distance=-0.18cm, label={\tiny $\Qt_Z^{(d')}$}]{4} 
\Edge[lw=1pt, position={left=1mm}](3)(4)

\Vertex[x=0,y=-3.5,opacity=0,position=below left,style={fill=blue,fill opacity=0.3},label=\tiny{\textcolor{blue}{$\widehat{\Wt}{}^{(d)}_{\text{D}}$}}, size=0.5, position= below right, distance=-0.2cm ]{5}
\Edge[lw=1pt, position={left=1mm}](4)(5)

\Vertex[x=0,y=-4.1,style={draw=none, fill=none}, size=.1]{W6}
\Edge[lw=1pt,position={left=1mm}](5)(W6)

\Vertex[x=0,y=-4.6,style={draw=none, fill=none}, size=.1]{W7}
\Edge[style={dotted}](W6)(W7)

\Vertex[x=0,y=-5.2, opacity=0,size=0.5]{8}
\Edge[lw=1pt, position={left=1mm}](W7)(8)

\Vertex[x=-0.75,y=-0,style={draw=none, fill=none}, size=.1]{I1}
\Vertex[x=-0.75,y=-1.7,style={draw=none, fill=none}, size=.1]{I3}
\Vertex[x=-0.75,y=-2.6,style={draw=none, fill=none}, size=.1]{I4}
\Vertex[x=-0.75,y=-3.5,style={draw=none, fill=none}, size=.1]{I5}
\Vertex[x=-0.75,y=-5.2,style={draw=none, fill=none}, size=.1]{I8}

\Edge[lw=1pt, position={left=1mm}](I1)(1)
\Edge[lw=1pt, position={left=1mm}](I3)(3)
\Edge[lw=1pt, position={left=1mm}](I4)(4)
\Edge[lw=1pt, position={left=1mm}](I5)(5)
\Edge[lw=1pt, position={left=1mm}](I8)(8)

\node[rectangle, rounded corners, draw=none, fill=none, fill opacity=0.3, text opacity = 1, below of=8, align=center, xshift=-0.5cm, yshift=0.5cm] {\tiny Thm. \ref{thm:VTN_obtain_WLSE}};
\node[rectangle, rounded corners, draw=none, fill=none, fill opacity=0.3, text opacity = 1, below of=8, align=center] (9) {\rotatebox{90}{$\,=$} \\ \tiny $\hat{\tilde{\w}}{}_{\text{LSE}}^{\text{TT}}$};

\end{scope}%
\begin{scope}[xshift=2.8cm, scale=0.85]
\Vertex[x=0,y=0,opacity=0, size=0.5]{1} 
\Vertex[x=0,y=-0.6,style={draw=none, fill=none}, size=.01]{W1}
\Edge[lw=1pt, position={left=1mm}](1)(W1)

\Vertex[x=0,y=-1.1,style={draw=none,fill=none}, size=.1]{W2}
\Edge[style={dotted}](W1)(W2)

\Vertex[x=0,y=-1.7,opacity=0, size=0.5]{3}
\Edge[lw=1pt, position={left=1mm}](W2)(3)

\Vertex[x=0,y=-2.6,opacity=0, size=0.5, position= below right, distance=-0.23cm,style={fill=blue,fill opacity=0.3}, label=\tiny{\textcolor{blue}{$\widehat{\Wt}{}^{(d')}_{\text{LSE}}$}}]{4} 
\Edge[lw=1pt, position={left=1mm}](3)(4)

\Vertex[x=0,y=-3.5,opacity=0, size=0.5]{5} 
\Edge[lw=1pt, position={left=1mm}](4)(5)

\Vertex[x=0,y=-4.1,style={draw=none, fill=none}, size=.1]{W6}
\Edge[lw=1pt,position={left=1mm}](5)(W6)

\Vertex[x=0,y=-4.6,style={draw=none, fill=none}, size=.1]{W7}
\Edge[style={dotted}](W6)(W7)

\Vertex[x=0,y=-5.2, opacity=0,size=0.5]{8}
\Edge[lw=1pt, position={left=1mm}](W7)(8)

\Vertex[x=-0.75,y=-0,style={draw=none, fill=none}, size=.1]{I1}
\Vertex[x=-0.75,y=-1.7,style={draw=none, fill=none}, size=.1]{I3}
\Vertex[x=-0.75,y=-2.6,style={draw=none, fill=none}, size=.1]{I4}
\Vertex[x=-0.75,y=-3.5,style={draw=none, fill=none}, size=.1]{I5}
\Vertex[x=-0.75,y=-5.2,style={draw=none, fill=none}, size=.1]{I8}

\Edge[lw=1pt, position={left=1mm}](I1)(1)
\Edge[lw=1pt, position={left=1mm}](I3)(3)
\Edge[lw=1pt, position={left=1mm}](I4)(4)
\Edge[lw=1pt, position={left=1mm}](I5)(5)
\Edge[lw=1pt, position={left=1mm}](I8)(8)
\node[rectangle, rounded corners, draw=none, fill=none, fill opacity=0.3, text opacity = 1, below of=8, align=center, xshift=-0.5cm, yshift=0.5cm] {\tiny Thm. \ref{thm:VTN_obtain_WLSE_2} };
\node[rectangle, rounded corners, draw=none, fill=none, fill opacity=0.3, text opacity = 1, below of=8, align=center] {\rotatebox{90}{$\,=$} \\ \tiny $\hat{\w}{}_{\text{LSE}}^{\text{TT}}$};
\end{scope}%
\begin{scope}[xshift=4.4cm, scale=0.85]
\Vertex[x=0,y=0,opacity=0, size=0.5]{1} 
\Vertex[x=0,y=-0.6,style={draw=none, fill=none}, size=.01]{W1}
\Edge[lw=1pt, position={left=1mm}](1)(W1)

\Vertex[x=0,y=-1.1,style={draw=none,fill=none}, size=.1]{W2}
\Edge[style={dotted}](W1)(W2)

\Vertex[x=0,y=-1.7,opacity=0, size=0.5]{3}
\Edge[lw=1pt, position={left=1mm}](W2)(3)

\Vertex[x=0,y=-2.6,opacity=0, size=0.5, style={fill=blue,fill opacity=0.3}, position= below right, distance=-0.23cm,label=\tiny{\textcolor{blue}{$\widehat{\Wt}{}^{(d')}_{D+1}$}}]{4}
\Edge[lw=1pt, position={left=1mm}](3)(4)

\Vertex[x=0,y=-3.5,opacity=0, size=0.5]{5} 
\Edge[lw=1pt, position={left=1mm}](4)(5)

\Vertex[x=0,y=-4.1,style={draw=none, fill=none}, size=.1]{W6}
\Edge[lw=1pt,position={left=1mm}](5)(W6)

\Vertex[x=0,y=-4.6,style={draw=none, fill=none}, size=.1]{W7}
\Edge[style={dotted}](W6)(W7)

\Vertex[x=0,y=-5.2, opacity=0,size=0.5]{8}
\Edge[lw=1pt, position={left=1mm}](W7)(8)

\Vertex[x=-0.75,y=-0,style={draw=none, fill=none}, size=.1]{I1}
\Vertex[x=-0.75,y=-1.7,style={draw=none, fill=none}, size=.1]{I3}
\Vertex[x=-0.75,y=-2.6,style={draw=none, fill=none}, size=.1]{I4}
\Vertex[x=-0.75,y=-3.5,style={draw=none, fill=none}, size=.1]{I5}
\Vertex[x=-0.75,y=-5.2,style={draw=none, fill=none}, size=.1]{I8}

\Edge[lw=1pt, position={left=1mm}](I1)(1)
\Edge[lw=1pt, position={left=1mm}](I3)(3)
\Edge[lw=1pt, position={left=1mm}](I4)(4)
\Edge[lw=1pt, position={left=1mm}](I5)(5)
\Edge[lw=1pt, position={left=1mm}](I8)(8)
\node[rectangle, rounded corners, draw=none, fill=none, fill opacity=0.3, text opacity = 1, below of=8, align=center] {\rotatebox{90}{$\,=$} \\ \tiny $\hat{\w}_{D+1}^{\text{TT}}$};
\end{scope}%
\begin{scope}[xshift=6cm, scale=0.85]
\Vertex[x=0,y=0,opacity=0, size=0.5]{1} 
\Vertex[x=0,y=-0.6,style={draw=none, fill=none}, size=.01]{W1}
\Edge[lw=1pt, position={left=1mm}](1)(W1)

\Vertex[x=0,y=-1.1,style={draw=none,fill=none}, size=.1]{W2}
\Edge[style={dotted}](W1)(W2)

\Vertex[x=0,y=-1.7,opacity=0, size=0.5]{3}
\Edge[lw=1pt, position={left=1mm}](W2)(3)

\Vertex[x=0,y=-2.6,opacity=0, size=0.5, style={fill=blue,fill opacity=0.3},position=right,distance=-0.1cm,label=\textcolor{blue}{\tiny$\widehat{\Wt}{}^{(d')}_{\Delta}$}]{4} 
\Vertex[x=0.3,y=-3,opacity=0, size=0.5, style={fill=none, draw=none},position=right,distance=-0.1cm,label=\textcolor{blue}{\tiny \rotatebox{90}{$=$}}]{X1}
\Vertex[x=0,y=-3.4,opacity=0, size=0.5, style={fill=none, draw=none},position=right,distance=-0.1cm,label=\textcolor{blue}{\tiny$\widehat{\Wt}{}^{(d')}_{D+1}$}]{X2} 
\Vertex[x=0.22,y=-3.8,opacity=0, size=0.5, style={fill=none, draw=none},position=right,distance=-0.1cm,label=\textcolor{blue}{\tiny \rotatebox{90}{$-$}}]{X3} 
\Vertex[x=0,y=-4.2,opacity=0, size=0.5,style={fill=none, draw=none},position=right,distance=-0.1cm, label=\textcolor{blue}{\tiny$\widehat{\Wt}{}^{(d')}_{\text{LSE}}$}]{X4} 

\Edge[lw=1pt, position={left=1mm}](3)(4)

\Vertex[x=0,y=-3.5,opacity=0, size=0.5]{5} 
\Edge[lw=1pt, position={left=1mm}](4)(5)

\Vertex[x=0,y=-4.1,style={draw=none, fill=none}, size=.1]{W6}
\Edge[lw=1pt,position={left=1mm}](5)(W6)

\Vertex[x=0,y=-4.6,style={draw=none, fill=none}, size=.1]{W7}
\Edge[style={dotted}](W6)(W7)

\Vertex[x=0,y=-5.2, opacity=0,size=0.5]{8}
\Edge[lw=1pt, position={left=1mm}](W7)(8)

\Vertex[x=-0.75,y=-0,style={draw=none, fill=none}, size=.1]{I1}
\Vertex[x=-0.75,y=-1.7,style={draw=none, fill=none}, size=.1]{I3}
\Vertex[x=-0.75,y=-2.6,style={draw=none, fill=none}, size=.1]{I4}
\Vertex[x=-0.75,y=-3.5,style={draw=none, fill=none}, size=.1]{I5}
\Vertex[x=-0.75,y=-5.2,style={draw=none, fill=none}, size=.1]{I8}

\Edge[lw=1pt, position={left=1mm}](I1)(1)
\Edge[lw=1pt, position={left=1mm}](I3)(3)
\Edge[lw=1pt, position={left=1mm}](I4)(4)
\Edge[lw=1pt, position={left=1mm}](I5)(5)
\Edge[lw=1pt, position={left=1mm}](I8)(8)
\node[rectangle, rounded corners, draw=none, fill=none, fill opacity=0.3, text opacity = 1, below of=8, align=center, xshift=-0.5cm, yshift=0.5cm] {\tiny Lem. \ref{lem:efficient_w_delta_TT}};
\node[rectangle, rounded corners, draw=none, fill=none, fill opacity=0.3, text opacity = 1, below of=8, align=center] {\rotatebox{90}{$\,=$} \\ \tiny $\hat{\w}_{\Delta}^{\text{TT}}$};
\end{scope}%
\begin{scope}[yshift=-0.7cm, xshift=-0.3cm, scale = 0.85]]
\node[draw=none] at (0.5,-1.8)  (1) {};
  \node[draw=none] at (1.5,-1.9) (2) [right of=1] {};
  \path[every node/.style={font=\sffamily\small},->,>=stealth',thick]
    (1) edge[bend left=60,color=green] node [above, xshift=0.18cm, yshift=-0.16cm] {\tiny Insert $\scriptscriptstyle \Qt_Z^{(d')}$} (2);
\end{scope}%
\begin{scope}[xshift=1.4cm, scale = 0.85]]

    \node[draw=none] at (0.5,-3.5)  (1) {};
     \node[draw=none] at (1.3,-2.6) (2) {};
    \path[every node/.style={font=\sffamily\small},->,>=stealth',thick] (1) edge[bend right=45, color=green] node[above,xshift=-0.25cm,yshift=-0.1cm] {\tiny LQ} (2);
\end{scope}%
\begin{scope}[xshift=2.65cm,yshift=-0.7cm,scale = 0.85]]
\node[draw=none] at (0.5,-1.8)  (1) {};
  \node[draw=none] at (1.5,-1.9) (2) [right of=1] {};
  \path[every node/.style={font=\sffamily\small},->,>=stealth',thick]
    (1) edge[bend left=60,color=green] node [above, yshift=-0.02cm] {\tiny Update } (2);
\end{scope}%
\begin{scope}[xshift=4.25cm,yshift=-0.7cm,scale = 0.85]]
\node[draw=none] at (0.5,-1.8)  (1) {};
  \node[draw=none] at (1.5,-1.9) (2) [right of=1] {};
  \path[every node/.style={font=\sffamily\small},->,>=stealth',thick]
    (1) edge[bend left=60,color=green] node [above, yshift=-0.02cm] {\tiny Substract } (2);
\end{scope}
\end{tikzpicture}
 \vspace{-0.3cm}
\caption{Increasing $D$ with weight-TT. \textbf{First TN:} $D$th order VTN in site-$d$-mixed canonical form. \textbf{Second TN:} Insertion of core $\scriptstyle \Qts^{(d')}$. \textbf{Third TN:} Shift norm to $d'$th core. \textbf{Fourth TN:} Update $d'$-th core. \textbf{Fifth TN:} Residual VTN model.}\label{fig:incD_TensorDiagram}
\end{figure}
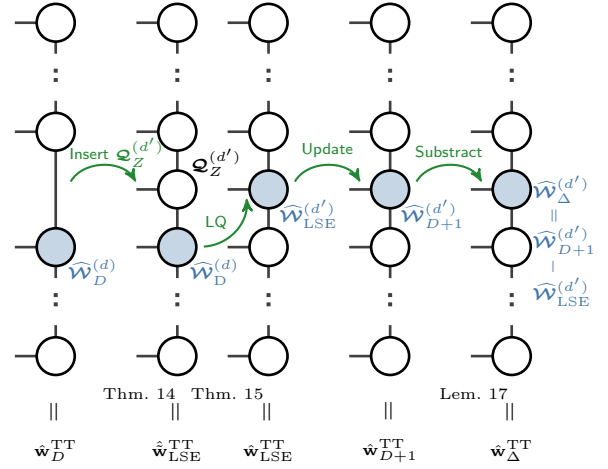
We assume the solution $\hat{\w}_D$ to the $D$th-order LS problem in Equation (\ref{eq:VM_LS_D}) to be given as the weight-TT $\hat{\w}_{D}^{\text{TT}}$ in Figure \ref{fig:incD_TensorDiagram}. In side-$d$-mixed canonical form, $\hat{\w}_{D}^{\text{TT}}$ is
\begin{align}
\hat{\w}{}_{D}^{\text{TT}} = \mathrm{TT}_{D}\left(\scriptstyle \hat{\Qt}^{(1)},\dots,\hat{\Qt}^{(d-1)},\hat{\Wt}^{(d)}_{D},\hat{\Qt}^{(d+1)},\dots,\hat{\Qt}^{(D)} \right). \label{eq:TT_level_0}
\end{align}
For the VTN formulation of the $D+1$th order LSE problem in Equation (\ref{eq:VM_LSE}), we construct the constraints just as in Section Section \ref{sec:LSE_VM} to zero out all additional parameters introduced by increasing the model order from $D$ to $D+1$.
Since computing the predictions $\hat{\y}$ in Equation (\ref{eq:MIMOVolterraAllN}) involves $D$ mode-$2$-products between a TT-core and the vector $\uv_n$ \cite{batselier2017tensor}, it is sufficient to impose the equality constraint $\C_{\text{TT}}$ on a single TT-core.
To facilitate establishing a connection between the $D$th-order and $D+1$th order VTN later on, we index the $D+1$th TT-cores of the weight-TT $\hat{\w}_{\text{LSE}}^{\text{TT}}$ as $\{,\dots,d-1,\ d',\ d,\dots,D\}$. 
Constraining the $d'$th core yields the ALS optimization step
\begin{align}
    \min_{\w_{\text{LSE}}^{(d')}} ||\y &- \U^{\otimes D+1}\hat{\Q}_{\setminus d'}\w_{\text{LSE}}^{(d')}|| \nonumber \\
    &\text{s.t.} \quad \C_{\text{TT}}\w_{\text{LSE}}^{(d')} = \mathbf{0}_{R_{d'}MR_{d'}}, \label{eq:VTN_LSE}
\end{align}
where $\w_{\text{LSE}}^{(d')} = \vect(\Wt_{\text{LSE}}^{(d')})$ with $\Wt_{\text{LSE}}^{(d')} \in \mathbb{R}^{R_{d'} \times I \times R_{d'}}$, and $\C_{\text{TT}} \in \mathbb{R}^{R_{d'}MR_{d'} \times R_{d'}IR_{d'}}$.
Similar to Section \ref{sec:LSE_VM}, we achieve the desired constraints by composing the rows $\C_{\text{TT}}$ of $R_{d'}MR_{d'}$ standard basis vectors $\ev_{R_{d'}IR_{d'},m}^\top$.
The reduced LSU problem in VTN format is given by
\begin{align}
\min_{\w_{\text{LSU}}^{(d')}} ||\y &- \U^{\otimes D+1}\hat{\Q}_{\setminus d'}\Z_{\text{TT}}\w_{\text{LSU}}^{(d')}||, \label{eq:VTN_LSU}
\end{align}
where the columns of $\Z_{\text{TT}} \in \mathbb{R}^{R_{d'}IR_{d'} \times R_{d'}^2}$ span $\mathcal{N}(\C_{\text{TT}})$; and $\w_{\text{LSU}}^{(d')} \in \mathbb{R}^{R_{d'}^2}$. We explicitly derive $\C_{\text{TT}}$ and $\Z_{\text{TT}}$ in Appendix \ref{app:CTTandZTT}.\\
\\
First, we transfer the results of Lemma \ref{lem:LSU_VM} to the VTN format.
Then, the remaining steps involve leverage the non-uniqueness and multilinearity of the TT decomposition. 
We design the TT-core $\Qt_\Z^{(d')}$, insert $\Qt_\Z^{(d')}$ into $\hat{\w}{}_{D}^{\text{TT}}$ in Equation (\ref{eq:TT_level_0}) and obtain $\hat{\widetilde{\w}}{}_{\text{LSE}}^{\text{TT}}$. 
Thus, $\Qt_\Z^{(d')}$ takes over the task of $\Z$ in Equation (\ref{eq:VM_LSU}) and is given by 
\begin{align}
   \Q_Z^{(d')}\,\big|_{[:,i,:]} := \begin{cases}
 \I_{R_{d'}}, \quad & \text{if $i = 1$},\\
 \mathbf{0}_{R_{d'} \times R_{d'}}, \quad & \text{if $2\leq i \leq I$}.
    \end{cases} \label{eq:incD_Zcore}
\end{align}
The construction of $\Qt_\Z^{(d')}$ is discussed in Appendix \ref{app:constructZTT}.
The resulting weight-TT $\hat{\widetilde{\w}}{}_{\text{LSE}}^{\text{TT}}$ is then given by  
 \begin{align}
    \hat{\widetilde{\w}}&{}_{\text{LSE}}^{\text{TT}} = \mathrm{TT}_{D+1}\left(\hat{\Qt}^{(1)},\dots\right. \nonumber \\
    & \dots,\left.\hat{\Qt}^{(d-1)},\Qt_\Z^{(d')},\hat{\Wt}^{(d)}_{D},\hat{\Qt}^{(d+1)},\dots,\hat{\Qt}^{(D)} \right). \label{eq:TT_level_1}
\end{align}
\begin{thm}\label{thm:VTN_obtain_WLSE}
    Inserting TT-core $\Qt_\Z^{(d')}$ defined in Equation (\ref{eq:incD_Zcore}) into the weight-TT $\hat{\w}{}_{D}^{\text{TT}}$ given in Equation (\ref{eq:TT_level_0}) is equivalent to computing $\hat{\widetilde{\w}}{}_{\text{LSE}}^{\text{TT}} = \Z\hat{\w}_D^{\text{TT}}$. 
\end{thm} 
Due to the non-uniqueness of the TT decomposition, we can obtain $\hat{\w}{}_{\text{LSE}}^{\text{TT}}$ by bringing $\hat{\widetilde{\w}}{}_{\text{LSE}}^{\text{TT}}$ into site-$d'$-mixed canonical form with the ALS-orthogonalization step.
Performing the LQ-decomposition of $\widehat{\W}{}^{(d)}_D|_\text{right}$ yields 
\begin{align}
\hat{\w}&{}_{\text{LSE}}^{\text{TT}} = \dots \nonumber \\
&\mathrm{TT}_{D+1}\left( \scriptstyle \hat{\Qt}^{(1)},\dots \hat{\Qt}^{(d-1)},\hat{\Wt}_{\text{LSE}}^{(d')},\hat{\Qt}^{(d)},\hat{\Qt}^{(d+1)},\dots,\hat{\Qt}^{(D)} \right), \label{eq:TT_level_2}
\end{align}
where $\hat{\Wt}_{\text{LSE}}^{(d')}$ and $\hat{\Qt}^{(d)}$ are the result of performing the ALS-orthogonalization step. 
The results of Corollary \ref{cor:identical_outputs} directly apply to the VTN.
\begin{thm}\label{thm:VTN_obtain_WLSE_2}
    Obtaining the TT-core $\hat{\Wt}_{\text{LSE}}^{(d')}$ in Equation 
    (\ref{eq:TT_level_2}) by bringing $\hat{\widetilde{\w}}{}_{\text{LSE}}^{\text{TT}}$ in Equation (\ref{eq:TT_level_1}) into site-$d'$-mixed canonical form is equivalent to solving the LSE problem defined in Equation (\ref{eq:VTN_LSE}).
\end{thm}
\begin{cor}\label{cor:identical_outputs}
The weight-TTs $\hat{\w}{}_{\text{LSE}}^{\text{TT}}$, $\hat{\tilde{\w}}{}_{\text{LSE}}^{\text{TT}}$ and $\hat{\w}{}_{\text{D}}^{\text{TT}}$ induce identical input-output behaviour $\hat{\y}_D = \U^{\otimes D}\hat{\w}{}_{D}^{\text{TT}}.$
\end{cor}
Lifting the constraints in Equation (\ref{eq:VTN_LSE}) results in the standard ALS-optimization step for the $d'$th TT-core of a $D+1$th order VTN, yielding the TT-core $\hat{\w}{}^{(d')}_{D+1}$ and the weight-TT $\hat{\w}{}_{\D+1}^{\text{TT}}$. 
Usually, computing $\hat{\w}_{\Delta}^{\text{TT}}$ with Equation (\ref{eq:updateVector_VM}) requires a rank-increasing procedure for addition, followed by a rounding step costing $\bigO(DIR^3)$ operations \cite{oseledets2011tensor}.
Since $\hat{\Q_{\setminus d'}}$ is identical for $\hat{\w}{}_{\text{LSE}}^{\text{TT}}$ based on $\hat{\w}{}_{\text{D+1}}^{\text{TT}}$, we can leverage the multilinearity of the TT-decomposition and compute $\hat{\w}_{\Delta}^{\text{TT}}$ with Lemma \ref{lem:efficient_w_delta_TT}.
\begin{lem}\label{lem:efficient_w_delta_TT}
    The weight-TT $\hat{\w}_{\Delta}^{\text{TT}}$ is given by 
    \begin{align}
        \hat{\w}_{\Delta}^{\text{TT}} = \hat{\Q_{\setminus d'}(\hat{\w}^{(d')}_{D+1} - \hat{\w}^{(d')}_{\text{LSE}})}\ \  .
    \end{align}
\end{lem}
The computations to increase $D$ to $D+1$ in VTN format cost $\bigO(M^3R^6)$. 
The derivation to increase $M$ in VTN format is largely similar.
The main difference is, that instead of adding one TT-core, we expand all TT-cores by the desired increase $M_\Delta$.
While increasing $D$ is performed with a single ASL-optimization step, increasing $M$ requires $D$ ALS-optimization steps and costs $\bigO(DM^3R^6)$.
The extension from MIMO to SISO does not influence our proposed method and is performed according to Remark \ref{rem:SISO2MIMO}.

\section{Deterministic Initialization and Hyperparameter Selection}
First, we demonstrate in Section \ref{sec:warmstart} the steps depicted in Figure \ref{fig:incD_TensorDiagram} can be employed as a deterministic initialization strategy for the SOTA VTN.
Then, we illustrate in Section \ref{sec:hyperparam_select1} and Section \ref{sec:hyperparam_select12} how to optimize the hyperparameters $D$ and $M$ using our proposed method.
\subsection{Deterministic Initialization Strategy}\label{sec:warmstart}
For the proposed deterministic initialization strategy, we start with a low-order deterministic model and then apply Figure \ref{fig:incD_TensorDiagram} iteratively until the desired order $D_\text{max}$ is reached. 
This approach introduces the hyperparameter \textit{sweeps}, which describes how many iterations (one TT-core update) are made between individual increases. 
It acts similar to the learning rate in gradient boosting \cite{friedman2001greedy}. 
We show in our experiments in Section \ref{sec:experiments}, how to tune $D$ and \textit{sweeps} simultaneously.
As confirmed by our experiments, the proposed deterministic initialization strategy has two advantages over the standard random initialization: we obtain fully deterministic results and improved performance, provided a good initial guess is used, due to the updates in conjugate directions. 
Since increasing the model order $D$ is achieved through a single core update, there are no additional computational costs compared to the SOTA VTN. \\
\\
We assume $D_\text{max}$, $M$ and $R$ are know. 
The first step is to compute a deterministic low-complexity model, e.g. of order $D=2$ or $D=1$, with the desired $M$ and $R$.
We use the algorithm presented in \cite{batselier2021enforcing}, which finds an optimal model in the least-squares sense. 
For small values of $D$, the algorithm scales linearly with $N$ with $\bigO(NI^2R_{\text{max}}^2)$, where $R_\text{max}$ denotes the symmetric upper bounds, which grow binomially with $D$ \cite[Theorem 2]{batselier2021enforcing}.
Note that choosing $D=1$ restricts the rank to $R=1$, which can be increased e.g. with the ALS-SD \cite{dolgov2014alternating} or the modified ALS \cite{batselier2017tensor}.
For $D=2$, we obtain the weight-TT 
\begin{align}
\mathrm{TT}_{2}\left( \scriptstyle \hat{\Wt}^{(1)}_{2},\hat{\Qt}^{(2)} \right), \label{eq:init_TT_0}
\end{align}
which serves as the deterministic initialization.
Next, we alternate between standard ALS iterations and one execution of Figure \ref{fig:incD_TensorDiagram} until the desired order $D_{\text{max}}$ is reached. 
The proposed method yields a VTN model of uniform rank $R$. 
To achieve e.g. a pyramidal rank structure, the modified ALS \cite{batselier2017tensor} can be used in place of the standard ALS iterations. 
\subsection{Efficient Grid-Search}\label{sec:hyperparam_select1} 
To perform an exhaustive grid search over an $M_\text{max} \times D_\text{set}$, we apply the deterministic initialization strategy outlined in Section \ref{sec:warmstart} for $M_\text{max} \times sweeps$ different models. 
This approach eliminates the required random restarts in SOTA VTN and reduces the grid size, since typically $D_\text{max} < sweeps$
We determine the optimal values $D_{\text{opt}}$, $M_{\text{opt}}$ and \textit{sweeps} on a validation data set and select the structure, that yields produces the lowest, transient free root mean square error (RMSE) given by
\begin{align}
    \text{RMSE} = \sqrt{\sum_{n=1}^N \frac{\left(y(n) - \hat{y}(n)\right)^2}{N}}.
\end{align}
Alternatively, we combine the RMSE with the idea behind the L-curve \cite{hansen1993use} and choose the choose the model structure as a trade-off between the RMSE on the validation data and the norm of the weight-TT.
\subsection{Automatic Hyperparameter Selection}\label{sec:hyperparam_select12} 
The incremental increase of $D$ and $M$ enables to simultaneous optimize both over a $M_\text{max} \times D_\text{set}$ in a single run with a fixed value for \textit{sweeps}. 
At each increase step, we increase both $D$ and $M$ and select the residual model $\hat{\w}_{\Delta}$ that yields the higher Variance Accounted For (VAF).
The orthogonality of $\hat{\y}_{D}$ and $\hat{\y}_\Delta$ ensures that the residual VAF does not contained variance explained by the previous model configuration.
The VAF is computed by 
\begin{align}
    \text{VAF} = 1 - \frac{\mathrm{var}\, (\y - \hat{\y})}{\mathrm{var}\, (\y)}.
\end{align}

\section{Experiments}\label{sec:experiments}
In the first experiment, we use the synthetical dataset proposed in \cite{batselier2021enforcing}. We showcase a grid search over $D$ and $M$ according to Section \ref{sec:hyperparam_select1} as well as the automatic hyperparameter selection proposed in Section \ref{sec:hyperparam_select12}.
In the second experiment, we use the nonlinear benchmark Cascaded Tanks \cite{schoukens2016cascaded}.
We demonstrate how to use to perform a grid search over $D \times sweeps$.\\
\\
All experiments were performed on a laptop running a 4-core Intel(R) Core(TM) i7-8665U CPUat 1.9 Ghz with 8GB RAM.
The code is written in MATLAB and will be made available after publication. 
\begin{figure}
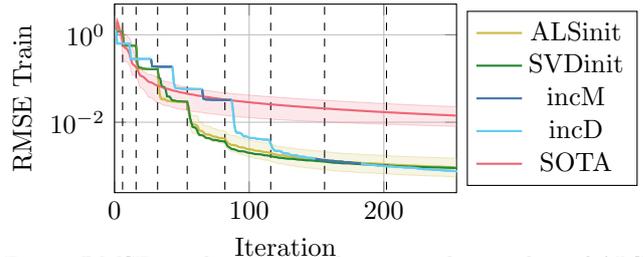

    \centering 

    \vspace{-0.5cm}
    \caption{RMSE on the training data over the number of ALS iterations. Comparing the SOTA VTN \cite{batselier2017tensor} (SOTA), increasing $D$ with deterministic (SVDinit) and random (ALSinit) initialization, simultaneously increasing $D$ (incD) and $M$ (incM) with deterministic initialization. The dashed line indicates increase steps for ALSinit and SVDinit. If a fixed value is applicable, $D = 7$, $M = 6$, $R = 3$ and $sweeps = 3$. }
    \label{fig:increaseD_nonoise_trainfill}
\end{figure}

\subsection{Synthetical Dataset}
We generate a $D = 7$, $M = 5$ SISO Volterra system based on \cite{batselier2021enforcing}.
We generate a $I \times 100$ matrix $\A$, whose exponentially decaying columns are given by 
\begin{align*}
    \A(:,i) = |\alpha_i| \exp(-\beta_i[1:I]),
\end{align*}
where $\alpha_i$, $\beta_i$ are sampled from a standard normal distribution and a uniform distribution in the range $[1,10]$, respectively.
We obtain the vectorized Volterra kernels $\w$ by computing the $D-1$-times column wise Kronecker product of $\A$ with itself, followed by summing over all columns.
We use Gaussian white noise as the input signal to generate $4500$ samples, splitting the data into $3000$ samples for training and $750$ each for validation and testing. We perform a grid search over $D = 2,\dots,10$ and $M =3,\dots,7$ with a fixed $R = 3$ and $sweeps = 3$.\\
\\
For the grid search according to Section \ref{sec:hyperparam_select1}, we compute $\mathrm{TT}_{2}$ in Equation (\ref{eq:init_TT_0}) using the approach in \cite{batselier2021enforcing} (SVDinit) and initialize it randomly with 10 random restarts (ALSinit).
This way, we disentangle the effects of the random initialization from those of the incremental updates, we compute for each $M$ $\mathrm{TT}_{2}$.
For the automatic hyperparameter determination (SVDinitDM) we compute $\mathrm{TT}_{2}$ in Equation (\ref{eq:init_TT_0}) using \cite{batselier2021enforcing} with $M = 3$. We compare the results to performing the grid search with the SOTA VTN \cite{batselier2017tensor}, stopping after 244 iterations with 10 random restarts for every grid point.\\
\\
Figure \ref{fig:increaseD_nonoise_trainfill} shows the RMSE on the training data for all four methods, fixing (if applicable) $D = 7$ and $M = 6$.
The results clearly indicate the primary factor in the performance improvement is the proposed update mechanism rather than the deterministic initialization. 
Moreover, the results demonstrate that conjugate update directions achieve a larger training error reduction than a standard ALS iteration. 
\begin{figure}
    \centering
\pgfplotsset{
  /pgfplots/boxplot legend/.style={
    legend image code/.code={
      \draw [|-|,##1] (0,-2mm) -- (0,2mm);
      \node[rectangle, minimum size=1.5mm, draw, fill=##1] at (0,0) {};
    }
  }
}
\begin{tikzpicture}
\begin{groupplot}[
  group style={group size=1 by 2, vertical sep=0.5cm},
  width=8cm, height=3cm,
  grid=major,
  legend to name=sharedlegend,
  legend columns=-1
]
\nextgroupplot[
    boxplot/draw direction=y,
xmin = 5.5, xmax = 18.75,
ymode = log,
max space between ticks=20
minor y tick style={black},    ylabel={$||\mathbf{w}^{\text{TT}}||_2$},
    xtick={2,4.25,6.5,8.75,11,13.25,15.5,17.75,20},
    xticklabels={2,3,4,5,6,7,8,9,10},
    enlarge x limits=0.01,
    grid=major,
    cycle list={{red}},
    legend pos=north west,
    boxplot legend 
]
\addplot+[
boxplot prepared={draw position=1.556,
    lower whisker=1.342e+00,
    lower quartile=1.346e+00,
    median=1.349e+00,
    upper quartile=1.357e+00,
    upper whisker=1.360e+00,
}, draw=yellow, fill=yellow!30] coordinates {};
\addlegendentry{ALSinit}
\addplot+[
forget plot,
boxplot prepared={draw position=3.806,
    lower whisker=7.356e-01,
    lower quartile=7.390e-01,
    median=7.426e-01,
    upper quartile=7.455e-01,
    upper whisker=7.539e-01,
}, draw=yellow, fill=yellow!30] coordinates {};
\addplot+[
forget plot,
boxplot prepared={draw position=6.056,
    lower whisker=4.148e-01,
    lower quartile=4.320e-01,
    median=4.527e-01,
    upper quartile=4.652e-01,
    upper whisker=4.793e-01,
}, draw=yellow, fill=yellow!30] coordinates {};
\addplot+[
forget plot,
only marks, mark=*, mark options={yellow}
] coordinates {(6.056, 5.920e-01)};
\addplot+[
forget plot,
boxplot prepared={draw position=8.306,
    lower whisker=2.855e-01,
    lower quartile=3.172e-01,
    median=3.603e-01,
    upper quartile=4.002e-01,
    upper whisker=4.425e-01,
}, draw=yellow, fill=yellow!30] coordinates {};
\addplot+[
forget plot,
boxplot prepared={draw position=10.556,
    lower whisker=2.483e-01,
    lower quartile=2.589e-01,
    median=2.837e-01,
    upper quartile=2.968e-01,
    upper whisker=3.163e-01,
}, draw=yellow, fill=yellow!30] coordinates {};
\addplot+[
forget plot,
boxplot prepared={draw position=12.806,
    lower whisker=2.599e-01,
    lower quartile=2.686e-01,
    median=2.853e-01,
    upper quartile=3.083e-01,
    upper whisker=3.245e-01,
}, draw=yellow, fill=yellow!30] coordinates {};
\addplot+[
forget plot,
boxplot prepared={draw position=15.056,
    lower whisker=2.982e-01,
    lower quartile=3.141e-01,
    median=3.419e-01,
    upper quartile=4.059e-01,
    upper whisker=5.011e-01,
}, draw=yellow, fill=yellow!30] coordinates {};
\addplot+[
forget plot,
boxplot prepared={draw position=17.306,
    lower whisker=3.616e-01,
    lower quartile=6.486e-01,
    median=8.181e-01,
    upper quartile=9.928e-01,
    upper whisker=1.124e+00,
}, draw=yellow, fill=yellow!30] coordinates {};
\addplot+[
forget plot,
boxplot prepared={draw position=19.556,
    lower whisker=8.088e-01,
    lower quartile=2.297e+00,
    median=2.851e+00,
    upper quartile=5.604e+00,
    upper whisker=1.540e+01,
}, draw=yellow, fill=yellow!30] coordinates {};
\addplot+[
boxplot prepared={draw position=2.444,
    lower whisker=1.382e+00,
    lower quartile=1.382e+00,
    median=1.382e+00,
    upper quartile=1.382e+00,
    upper whisker=1.383e+00,
}, draw=red, fill=red!30] coordinates {};
\addlegendentry{SOTA}
\addplot+[
forget plot,
only marks, mark=*, mark options={red}
] coordinates {(2.444, 1.379e+00)};
\addplot+[
forget plot,
only marks, mark=*, mark options={red}
] coordinates {(2.444, 1.378e+00)};
\addplot+[
forget plot,
boxplot prepared={draw position=4.694,
    lower whisker=7.923e-01,
    lower quartile=7.957e-01,
    median=8.286e-01,
    upper quartile=8.556e-01,
    upper whisker=9.150e-01,
}, draw=red, fill=red!30] coordinates {};
\addplot+[
forget plot,
boxplot prepared={draw position=6.944,
    lower whisker=4.904e-01,
    lower quartile=5.265e-01,
    median=5.662e-01,
    upper quartile=6.677e-01,
    upper whisker=7.335e-01,
}, draw=red, fill=red!30] coordinates {};
\addplot+[
forget plot,
boxplot prepared={draw position=9.194,
    lower whisker=2.657e-01,
    lower quartile=2.827e-01,
    median=3.924e-01,
    upper quartile=4.535e-01,
    upper whisker=6.337e-01,
}, draw=red, fill=red!30] coordinates {};
\addplot+[
forget plot,
boxplot prepared={draw position=11.444,
    lower whisker=2.317e-01,
    lower quartile=2.448e-01,
    median=2.702e-01,
    upper quartile=3.127e-01,
    upper whisker=3.816e-01,
}, draw=red, fill=red!30] coordinates {};
\addplot+[
forget plot,
boxplot prepared={draw position=13.694,
    lower whisker=2.839e-01,
    lower quartile=3.049e-01,
    median=3.165e-01,
    upper quartile=3.503e-01,
    upper whisker=3.775e-01,
}, draw=red, fill=red!30] coordinates {};
\addplot+[
forget plot,
boxplot prepared={draw position=15.944,
    lower whisker=2.842e-01,
    lower quartile=3.098e-01,
    median=3.588e-01,
    upper quartile=4.177e-01,
    upper whisker=4.946e-01,
}, draw=red, fill=red!30] coordinates {};
\addplot+[
forget plot,
boxplot prepared={draw position=18.194,
    lower whisker=2.905e-01,
    lower quartile=3.464e-01,
    median=3.750e-01,
    upper quartile=4.170e-01,
    upper whisker=4.373e-01,
}, draw=red, fill=red!30] coordinates {};
\addplot+[
forget plot,
boxplot prepared={draw position=20.444,
    lower whisker=4.158e-01,
    lower quartile=4.470e-01,
    median=5.748e-01,
    upper quartile=6.382e-01,
    upper whisker=9.193e-01,
}, draw=red, fill=red!30] coordinates {};
\addplot+[
boxplot prepared={draw position=2.000,
    lower whisker=1.322e+00,
    lower quartile=1.322e+00,
    median=1.322e+00,
    upper quartile=1.322e+00,
    upper whisker=1.322e+00,
}, draw=green, fill=green!30] coordinates {};
\addplot+[
forget plot,
boxplot prepared={draw position=4.250,
    lower whisker=2.250e+00,
    lower quartile=2.250e+00,
    median=2.250e+00,
    upper quartile=2.250e+00,
    upper whisker=2.250e+00,
}, draw=green, fill=green!30] coordinates {};
\addplot+[
forget plot,
boxplot prepared={draw position=6.500,
    lower whisker=1.093e+00,
    lower quartile=1.093e+00,
    median=1.093e+00,
    upper quartile=1.093e+00,
    upper whisker=1.093e+00,
}, draw=green, fill=green!30] coordinates {};
\addplot+[
forget plot,
boxplot prepared={draw position=8.750,
    lower whisker=4.792e-01,
    lower quartile=4.792e-01,
    median=4.792e-01,
    upper quartile=4.792e-01,
    upper whisker=4.792e-01,
}, draw=green, fill=green!30] coordinates {};
\addplot+[
forget plot,
boxplot prepared={draw position=11.000,
    lower whisker=3.027e-01,
    lower quartile=3.027e-01,
    median=3.027e-01,
    upper quartile=3.027e-01,
    upper whisker=3.027e-01,
}, draw=green, fill=green!30] coordinates {};
\addplot+[
forget plot,
boxplot prepared={draw position=13.250,
    lower whisker=3.501e-01,
    lower quartile=3.501e-01,
    median=3.501e-01,
    upper quartile=3.501e-01,
    upper whisker=3.501e-01,
}, draw=green, fill=green!30] coordinates {};
\addplot+[
forget plot,
boxplot prepared={draw position=15.500,
    lower whisker=4.835e-01,
    lower quartile=4.835e-01,
    median=4.835e-01,
    upper quartile=4.835e-01,
    upper whisker=4.835e-01,
}, draw=green, fill=green!30] coordinates {};
\addplot+[
forget plot,
boxplot prepared={draw position=17.750,
    lower whisker=1.792e+00,
    lower quartile=1.792e+00,
    median=1.792e+00,
    upper quartile=1.792e+00,
    upper whisker=1.792e+00,
}, draw=green, fill=green!30] coordinates {};
\addplot+[
forget plot,
boxplot prepared={draw position=20.000,
    lower whisker=2.856e+01,
    lower quartile=2.856e+01,
    median=2.856e+01,
    upper quartile=2.856e+01,
    upper whisker=2.856e+01,
}, draw=green, fill=green!30] coordinates {};
\addlegendimage{boxplot legend, green}
\addlegendentry{SVDinit}
\nextgroupplot[
    boxplot/draw direction=y,
xmin = 5.5, xmax = 18.75,
ymode = log,
max space between ticks=20
minor y tick style={black},    ylabel={RMSE Val},
    xtick={2,4.25,6.5,8.75,11,13.25,15.5,17.75,20},
    xticklabels={2,3,4,5,6,7,8,9,10},
    enlarge x limits=0.01,
    grid=major,
    cycle list={{red}},
    legend pos=south west,
    boxplot legend 
]
\addplot+[
boxplot prepared={draw position=1.556,
    lower whisker=8.719e-01,
    lower quartile=8.731e-01,
    median=8.735e-01,
    upper quartile=8.756e-01,
    upper whisker=8.773e-01,
}, draw=yellow, fill=yellow!30] coordinates {};
\addlegendentry{ALSinit}
\addplot+[
forget plot,
boxplot prepared={draw position=3.806,
    lower whisker=4.558e-01,
    lower quartile=4.584e-01,
    median=4.592e-01,
    upper quartile=4.610e-01,
    upper whisker=4.614e-01,
}, draw=yellow, fill=yellow!30] coordinates {};
\addplot+[
forget plot,
boxplot prepared={draw position=6.056,
    lower whisker=1.656e-01,
    lower quartile=1.671e-01,
    median=1.684e-01,
    upper quartile=1.700e-01,
    upper whisker=1.720e-01,
}, draw=yellow, fill=yellow!30] coordinates {};
\addplot+[
forget plot,
boxplot prepared={draw position=8.306,
    lower whisker=3.594e-02,
    lower quartile=3.694e-02,
    median=3.784e-02,
    upper quartile=3.940e-02,
    upper whisker=4.132e-02,
}, draw=yellow, fill=yellow!30] coordinates {};
\addplot+[
forget plot,
boxplot prepared={draw position=10.556,
    lower whisker=4.141e-03,
    lower quartile=5.329e-03,
    median=5.950e-03,
    upper quartile=7.134e-03,
    upper whisker=1.020e-02,
}, draw=yellow, fill=yellow!30] coordinates {};
\addplot+[
forget plot,
boxplot prepared={draw position=12.806,
    lower whisker=1.281e-03,
    lower quartile=1.735e-03,
    median=2.642e-03,
    upper quartile=3.389e-03,
    upper whisker=4.754e-03,
}, draw=yellow, fill=yellow!30] coordinates {};
\addplot+[
forget plot,
boxplot prepared={draw position=15.056,
    lower whisker=9.418e-04,
    lower quartile=1.359e-03,
    median=1.554e-03,
    upper quartile=2.136e-03,
    upper whisker=3.759e-03,
}, draw=yellow, fill=yellow!30] coordinates {};
\addplot+[
forget plot,
boxplot prepared={draw position=17.306,
    lower whisker=8.754e-04,
    lower quartile=1.124e-03,
    median=1.351e-03,
    upper quartile=1.861e-03,
    upper whisker=2.211e-03,
}, draw=yellow, fill=yellow!30] coordinates {};
\addplot+[
forget plot,
boxplot prepared={draw position=19.556,
    lower whisker=9.084e-04,
    lower quartile=1.041e-03,
    median=1.105e-03,
    upper quartile=1.606e-03,
    upper whisker=2.739e-03,
}, draw=yellow, fill=yellow!30] coordinates {};
\addplot+[
boxplot prepared={draw position=2.000,
    lower whisker=8.737e-01,
    lower quartile=8.737e-01,
    median=8.737e-01,
    upper quartile=8.737e-01,
    upper whisker=8.737e-01,
}, draw=green, fill=green!30] coordinates {};
\addlegendentry{SVDinit}
\addplot+[
forget plot,
boxplot prepared={draw position=4.250,
    lower whisker=4.575e-01,
    lower quartile=4.575e-01,
    median=4.575e-01,
    upper quartile=4.575e-01,
    upper whisker=4.575e-01,
}, draw=green, fill=green!30] coordinates {};
\addplot+[
forget plot,
boxplot prepared={draw position=6.500,
    lower whisker=1.687e-01,
    lower quartile=1.687e-01,
    median=1.687e-01,
    upper quartile=1.687e-01,
    upper whisker=1.687e-01,
}, draw=green, fill=green!30] coordinates {};
\addplot+[
forget plot,
boxplot prepared={draw position=8.750,
    lower whisker=4.106e-02,
    lower quartile=4.106e-02,
    median=4.106e-02,
    upper quartile=4.106e-02,
    upper whisker=4.106e-02,
}, draw=green, fill=green!30] coordinates {};
\addplot+[
forget plot,
boxplot prepared={draw position=11.000,
    lower whisker=5.712e-03,
    lower quartile=5.712e-03,
    median=5.712e-03,
    upper quartile=5.712e-03,
    upper whisker=5.712e-03,
}, draw=green, fill=green!30] coordinates {};
\addplot+[
forget plot,
boxplot prepared={draw position=13.250,
    lower whisker=1.876e-03,
    lower quartile=1.876e-03,
    median=1.876e-03,
    upper quartile=1.876e-03,
    upper whisker=1.876e-03,
}, draw=green, fill=green!30] coordinates {};
\addplot+[
forget plot,
boxplot prepared={draw position=15.500,
    lower whisker=1.404e-03,
    lower quartile=1.404e-03,
    median=1.404e-03,
    upper quartile=1.404e-03,
    upper whisker=1.404e-03,
}, draw=green, fill=green!30] coordinates {};
\addplot+[
forget plot,
boxplot prepared={draw position=17.750,
    lower whisker=1.183e-03,
    lower quartile=1.183e-03,
    median=1.183e-03,
    upper quartile=1.183e-03,
    upper whisker=1.183e-03,
}, draw=green, fill=green!30] coordinates {};
\addplot+[
forget plot,
boxplot prepared={draw position=20.000,
    lower whisker=1.115e-03,
    lower quartile=1.115e-03,
    median=1.115e-03,
    upper quartile=1.115e-03,
    upper whisker=1.115e-03,
}, draw=green, fill=green!30] coordinates {};
\addplot+[
boxplot prepared={draw position=2.444,
    lower whisker=8.734e-01,
    lower quartile=8.734e-01,
    median=8.734e-01,
    upper quartile=8.735e-01,
    upper whisker=8.735e-01,
}, draw=red, fill=red!30] coordinates {};
\addlegendentry{SOTA}
\addplot+[
forget plot,
only marks, mark=*, mark options={red}
] coordinates {(2.444, 8.737e-01)};
\addplot+[
forget plot,
only marks, mark=*, mark options={red}
] coordinates {(2.444, 8.738e-01)};
\addplot+[
forget plot,
boxplot prepared={draw position=4.694,
    lower whisker=4.585e-01,
    lower quartile=4.603e-01,
    median=4.607e-01,
    upper quartile=4.616e-01,
    upper whisker=4.628e-01,
}, draw=red, fill=red!30] coordinates {};
\addplot+[
forget plot,
boxplot prepared={draw position=6.944,
    lower whisker=1.664e-01,
    lower quartile=1.675e-01,
    median=1.682e-01,
    upper quartile=1.689e-01,
    upper whisker=1.704e-01,
}, draw=red, fill=red!30] coordinates {};
\addplot+[
forget plot,
boxplot prepared={draw position=9.194,
    lower whisker=3.790e-02,
    lower quartile=3.893e-02,
    median=3.980e-02,
    upper quartile=4.027e-02,
    upper whisker=4.209e-02,
}, draw=red, fill=red!30] coordinates {};
\addplot+[
forget plot,
boxplot prepared={draw position=11.444,
    lower whisker=6.023e-03,
    lower quartile=7.411e-03,
    median=9.684e-03,
    upper quartile=1.236e-02,
    upper whisker=1.544e-02,
}, draw=red, fill=red!30] coordinates {};
\addplot+[
forget plot,
boxplot prepared={draw position=13.694,
    lower whisker=1.859e-02,
    lower quartile=2.107e-02,
    median=2.637e-02,
    upper quartile=3.287e-02,
    upper whisker=3.466e-02,
}, draw=red, fill=red!30] coordinates {};
\addplot+[
forget plot,
boxplot prepared={draw position=15.944,
    lower whisker=2.488e-02,
    lower quartile=4.587e-02,
    median=5.641e-02,
    upper quartile=7.142e-02,
    upper whisker=1.322e-01,
}, draw=red, fill=red!30] coordinates {};
\addplot+[
forget plot,
only marks, mark=*, mark options={red}
] coordinates {(15.944, 1.487e-01)};
\addplot+[
forget plot,
boxplot prepared={draw position=18.194,
    lower whisker=6.124e-02,
    lower quartile=6.804e-02,
    median=8.088e-02,
    upper quartile=9.524e-02,
    upper whisker=1.256e-01,
}, draw=red, fill=red!30] coordinates {};
\addplot+[
forget plot,
boxplot prepared={draw position=20.444,
    lower whisker=1.125e-01,
    lower quartile=1.512e-01,
    median=2.617e-01,
    upper quartile=4.089e-01,
    upper whisker=6.013e-01,
}, draw=red, fill=red!30] coordinates {};
\addlegendimage{boxplot legend, red}
\end{groupplot}
\node[anchor=south] at ($(group c1r1.north) + (0,0.1cm)$) {\ref{sharedlegend}};
\node[anchor=north] at ($(group c1r2.south) - (0,0.5cm)$) {$D$};
\end{tikzpicture}

    \caption{Boxplots in nonlinear log-scale \cite{cox2009speaking}; whiskers and outliers using a $1.5 \times$ interquantile range of logarithms. \textbf{Top:} Variance of the norm $||\hat{\w}{}^{\text{TT}}||_2$ over $D$. \textbf{Bottom:} Variance of the validation RMSE over $D$. $M = 6$, $sweeps = 3$ and $R = 3.$}.
    \label{fig:increaseD_noise_boxplot}
\end{figure}
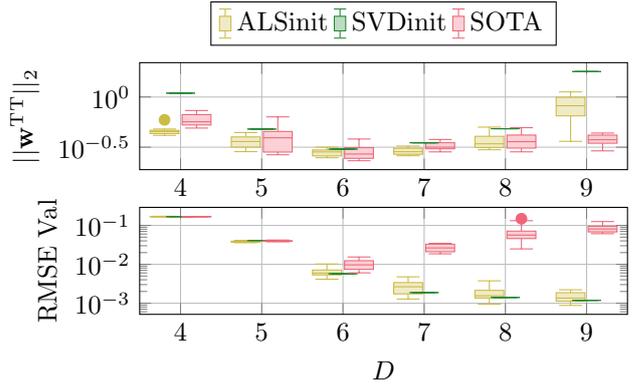
Instead of using the standard L-curve \cite{hansen1993use}, we plot the variance of the norm $||\hat{\w}{}^{\text{TT}}||_2$ and the RMSE on the validation data in Figure \ref{fig:increaseD_noise_boxplot}.
For ALSinit and SVDinit norm increases significantly as soon as the true underlying model order of $D=7$ is exceeded. 
This results from the non-zero parameters, that are added in the first ALS-optimization step after the increase step.
They compensate for the additional directions in $\hat{\Qt}{}^{(d')}$, which is obtained from the LQ-decomposition of $\widehat{\W}{}^{(d')}_{D+1}|_\text{right}$.
\begin{table}[]
    \centering
\begin{tabular}{lcccccc}
\toprule
\multirow{2}{*}{Method} & Selected & \multicolumn{2}{c}{Runtime (s)} & \multicolumn{2}{c}{RMSE (e-3)}\\
& D, M & Val. & Test & Val. & Test \\
\midrule
ALSinit & 7, 5 & 222 & 0.74 & 1.9 & 3.4 \\
SVDinit & 7, 6 & 23 & 0.07 & \textbf{1.1} &\textbf{1.8} \\
SVDinitDM  & 7, 6 & 4 & 0.06 & 1.5 & 5.6 \\
SOTA VTN & 6, 5 & 1594 & 0.84 & 10.1 & 15.2 \\
\bottomrule
\end{tabular}
\vspace{0.25cm}
    \caption{Performing grid search over $D = 2,\dots,10$, $M = 3,\dots,7$.}
    \label{tab:synthetic_comparison}
\end{table}\\
\\
For SOTA VTN and SVDinitDM, we select the model order based on the validation RMSE. For ALSinit and SVDinit, we select it based on the trade-off between model complexity and validation RMSE in Figure \ref{fig:increaseD_noise_boxplot}. Table \ref{tab:synthetic_comparison} compares the runtime for the grid search, as well as the validation and test error for the selected model structure.
Compared to the SOTA VTN, the combined method SVDinitDM achieves the most significant speedup of almost $400\times$ with an $8 \times$ improvement in RMSE. 
The ALSinit yields a $7 \times$ speedup with a more $4 \times$ improvement in accuracy. 
The SVDinit obtains $70 \times$ speedup with the highest RMSE improvement of more than  $8 \times$ compared to SOTA VTN. 

\subsection{Cascaded Tanks Benchmark}
We showcase the selection of the model order $D$ and the hyperparameter \textit{sweeps} on the nonlinear system identification benchmark Cascaded Tanks \cite{schoukens2016cascaded}.
The dataset consists of a training and test set, with $N=1024$ data points each and a Signal-to-Noise Ratio (SNR) of around $40\,\text{dB}$.
This benchmark is commonly used in the Volterra series literature to showcase different regularization techniques \cite{stoddard2018regularized, birpoutsoukis2018efficient, dalla2021kernel}.
Unfortunately, these approaches provide only limited comparability to our proposed method, since they either assume $D$ and $M$ to be known \cite{dalla2021kernel} or present three models of different orders without performing previous model selection \cite{stoddard2018regularized, birpoutsoukis2018efficient}.\\
\\
We split the given training data in a $2:1$ split in training and validation data, and select $D$ and \textit{sweeps} based on the lowest validation RMSE.
Then, we retrain a model in the chosen structure on the full training data and evaluate the selected model structure on the transient-free part of the test data.
We set $M = 95$ and $R=1$.
\begin{figure}
\centering
   \begin{subfigure}[t]{0.225\textwidth}
        \centering
\begin{tikzpicture}
    \begin{axis}[
        width=4cm,
        height=4cm,
        colorbar,
        xlabel = {$D$}, ylabel = {sweeps},
        colorbar horizontal,
        colorbar style={at={(0.5,1.1)},anchor=south,ticklabel pos=upper, height=0.25cm,},
        colormap name = viridis,
    ]
    \addplot[matrix plot, point meta=explicit]
        coordinates {
(2,1) [1.3836] (3,1) [1.17008] (4,1) [1.10174] (5,1) [1.04095] (6,1) [1.03215] (7,1) [1.02653] (8,1) [1.02224] (9,1) [1.01811] (10,1) [1.0134] (11,1) [1.00759] (12,1) [1.00019] (13,1) [0.990846] (14,1) [0.979409] (15,1) [0.96616]

(2,2) [1.3836] (3,2) [1.09012] (4,2) [1.01365] (5,2) [0.97764] (6,2) [0.95566] (7,2) [0.941262] (8,2) [0.929546] (9,2) [0.920852] (10,2) [0.915725] (11,2) [0.914352] (12,2) [0.917617] (13,2) [0.928532] (14,2) [0.951704] (15,2) [0.987828]

(2,3) [1.29586] (3,3) [1.08594] (4,3) [1.00614] (5,3) [0.990547] (6,3) [0.994104] (7,3) [1.04622] (8,3) [1.15089] (9,3) [1.26842] (10,3) [1.36677] (11,3) [1.45224] (12,3) [1.5359] (13,3) [1.62195] (14,3) [1.70989] (15,3) [1.79774]

(2,4) [1.28833] (3,4) [1.0724] (4,4) [1.0884] (5,4) [1.16249] (6,4) [1.25533] (7,4) [1.31917] (8,4) [1.33523] (9,4) [1.38808] (10,4) [1.46791] (11,4) [1.55909] (12,4) [1.65448] (13,4) [1.74228] (14,4) [1.81508] (15,4) [1.8728]

(2,5) [1.26796] (3,5) [1.08024] (4,5) [1.12946] (5,5) [1.13315] (6,5) [1.19181] (7,5) [1.28774] (8,5) [1.38114] (9,5) [1.46398] (10,5) [1.54685] (11,5) [1.63001] (12,5) [1.70654] (13,5) [1.76988] (14,5) [1.81877] (15,5) [1.85582]

(2,6) [1.17569] (3,6) [1.03249] (4,6) [1.10647] (5,6) [1.19341] (6,6) [1.28086] (7,6) [1.37098] (8,6) [1.46891] (9,6) [1.64324] (10,6) [1.82147] (11,6) [1.98865] (12,6) [2.11861] (13,6) [2.2092] (14,6) [2.27339] (15,6) [2.32123]

(2,7) [1.18443] (3,7) [1.09674] (4,7) [1.22335] (5,7) [1.29645] (6,7) [1.33194] (7,7) [1.40842] (8,7) [1.58272] (9,7) [1.76961] (10,7) [1.9139] (11,7) [2.0146] (12,7) [2.09463] (13,7) [2.16419] (14,7) [2.22495] (15,7) [2.2757]

(2,8) [1.18655] (3,8) [1.02735] (4,8) [1.04001] (5,8) [1.06702] (6,8) [1.13941] (7,8) [1.22516] (8,8) [1.44227] (9,8) [1.66062] (10,8) [1.80368] (11,8) [1.89041] (12,8) [1.94727] (13,8) [1.98868] (14,8) [2.02137] (15,8) [2.04853]

(2,9) [1.18377] (3,9) [1.1269] (4,9) [1.23824] (5,9) [1.29822] (6,9) [1.30453] (7,9) [1.40288] (8,9) [1.55253] (9,9) [1.70435] (10,9) [1.82014] (11,9) [1.90643] (12,9) [1.97635] (13,9) [2.03718] (14,9) [2.09228] (15,9) [2.14284]

(2,10) [1.15469] (3,10) [1.03461] (4,10) [1.02746] (5,10) [1.10947] (6,10) [1.18673] (7,10) [1.42423] (8,10) [1.67081] (9,10) [1.82617] (10,10) [1.91432] (11,10) [1.9762] (12,10) [2.02547] (13,10) [2.06596] (14,10) [2.09934] (15,10) [2.12708]

        };
    \end{axis}
\end{tikzpicture}

        \vspace{-0.7cm}
        \caption{Validation RMSE.}
        \label{fig:gridDS_heatmap_val}
    \end{subfigure}%
    ~ 
    \begin{subfigure}[t]{0.225\textwidth}
        \centering
\begin{tikzpicture}
    \begin{axis}[
        width=4cm,
        height=4cm,
        colorbar,
        xlabel = {$D$}, ylabel = {sweeps},
        colorbar horizontal,
        colorbar style={at={(0.5,1.1)},anchor=south,ticklabel pos=upper,height=0.25cm,},
        colormap name = viridis,
    ]
    \addplot[matrix plot, point meta=explicit]
        coordinates {
(2,1) [0.812137] (3,1) [0.780539] (4,1) [0.698864] (5,1) [0.691882] (6,1) [0.685942] (7,1) [0.681815] (8,1) [0.67924] (9,1) [0.677076] (10,1) [0.675276] (11,1) [0.663121] (12,1) [0.66117] (13,1) [0.658731] (14,1) [0.656425] (15,1) [0.654158]

(2,2) [0.812137] (3,2) [0.735876] (4,2) [0.681523] (5,2) [0.653296] (6,2) [0.64583] (7,2) [0.639563] (8,2) [0.631931] (9,2) [0.623802] (10,2) [0.615894] (11,2) [0.610313] (12,2) [0.609387] (13,2) [0.613648] (14,2) [0.62152] (15,2) [0.630345]

(2,3) [0.76943] (3,3) [0.685538] (4,3) [0.65651] (5,3) [0.654295] (6,3) [0.66004] (7,3) [0.667769] (8,3) [0.673029] (9,3) [0.674946] (10,3) [0.674063] (11,3) [0.670763] (12,3) [0.665318] (13,3) [0.658964] (14,3) [0.653987] (15,3) [0.65151]

(2,4) [0.727359] (3,4) [0.665589] (4,4) [0.649401] (5,4) [0.657134] (6,4) [0.664943] (7,4) [0.668068] (8,4) [0.665392] (9,4) [0.657767] (10,4) [0.647939] (11,4) [0.641582] (12,4) [0.639811] (13,4) [0.639691] (14,4) [0.639646] (15,4) [0.639424]

(2,5) [0.697284] (3,5) [0.654974] (4,5) [0.649722] (5,5) [0.662813] (6,5) [0.66699] (7,5) [0.666716] (8,5) [0.661812] (9,5) [0.65133] (10,5) [0.641298] (11,5) [0.637581] (12,5) [0.636849] (13,5) [0.636669] (14,5) [0.636647] (15,5) [0.636765]

(2,6) [0.666292] (3,6) [0.644517] (4,6) [0.644022] (5,6) [0.659524] (6,6) [0.666168] (7,6) [0.668292] (8,6) [0.665994] (9,6) [0.657817] (10,6) [0.645702] (11,6) [0.639439] (12,6) [0.638086] (13,6) [0.637794] (14,6) [0.637789] (15,6) [0.63797]

(2,7) [0.660737] (3,7) [0.6439] (4,7) [0.648512] (5,7) [0.663573] (6,7) [0.667323] (7,7) [0.66585] (8,7) [0.655743] (9,7) [0.642525] (10,7) [0.637709] (11,7) [0.63687] (12,7) [0.636733] (13,7) [0.636895] (14,7) [0.637243] (15,7) [0.637702]

(2,8) [0.65703] (3,8) [0.644655] (4,8) [0.652579] (5,8) [0.666596] (6,8) [0.666808] (7,8) [0.659304] (8,8) [0.643493] (9,8) [0.63732] (10,8) [0.63635] (11,8) [0.636252] (12,8) [0.63651] (13,8) [0.63696] (14,8) [0.637511] (15,8) [0.638125]

(2,9) [0.654411] (3,9) [0.64602] (4,9) [0.65642] (5,9) [0.668779] (6,9) [0.664463] (7,9) [0.649951] (8,9) [0.638167] (9,9) [0.636313] (10,9) [0.636097] (11,9) [0.636367] (12,9) [0.636869] (13,9) [0.637477] (14,9) [0.638143] (15,9) [0.638865]

(2,10) [0.652724] (3,10) [0.647658] (4,10) [0.660061] (5,10) [0.670248] (6,10) [0.660346] (7,10) [0.642689] (8,10) [0.636865] (9,10) [0.636191] (10,10) [0.636356] (11,10) [0.636858] (12,10) [0.637494] (13,10) [0.638187] (14,10) [0.638931] (15,10) [0.639752]

        };
    \end{axis}
\end{tikzpicture}

        \vspace{-0.7cm}
        \caption{Test RMSE.}
        \label{fig:gridDS_heatmap_test}
    \end{subfigure}
    \caption{Color represents RMSE on validation and test data.}
    \label{fig:gridDS_heatmap}
\end{figure}
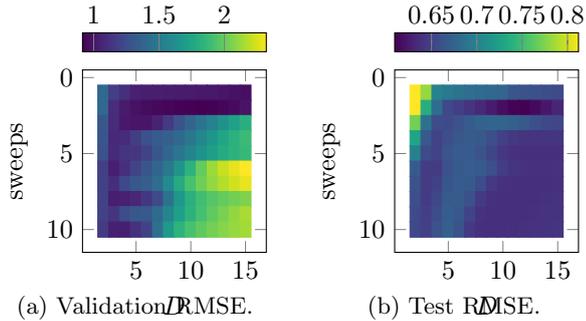
\\
\\
In Figure \ref{fig:gridDS_heatmap}, we analyze the trade-off between $D$ and \textit{sweeps} on both test and validation data.
We choose the deterministic initialization method and perform a grid search over $D = 2:15$ and $sweeps = 1:10$.
Setting small values for \textit{sweeps} results in a good model performance on unseen data, both in under- and overdramatized setting, indicating a regularizing effect of \textit{sweeps} similar to the learning rate in gradient boosting.
The results show a strong correlation between the validation and test error, such that we can choose the combination of $D$ and \textit{sweeps} based on the RMSE of the validation error.
Given that the interpolation threshold for the validation data is at $D=8$, while it is at $D=11$ for the test data, the clear difference between the validation and test error for larger values of $D$ and \textit{sweeps} is to be expected. 
The lowest validation RMSE is obtained for $D = 11$ and $sweeps = 2$ with $\mathrm{RMSE}_\text{val} = 0.91$. 
On the test data, this configuration yields $\mathrm{RMSE}_\text{test} = 0.61$.
The training and simulation time required to create Figure \ref{fig:gridDS_heatmap} is $74\,\text{s}$ for the validation and  $110\, \text{s} $ for the test data.
\\
\\
Next, we compare the proposed model order selection SVDinit procedure to the SOTA VTN implementation \cite{batselier2017tensor}, as well as the results reported  for the regularized Volterra series (RVS) in \cite{birpoutsoukis2018efficient} and the Laguerre basis functions (LBF)\cite{stoddard2018regularized}.
We design the search space for the SOTA VTN and SVDinit based on the respective methods' known properties. 
Overparametrized SOTA VTN models perform poorly on unseen data \cite{memmel2023bayesian}, so we restrict the search space to $D = 2,\dots,7$ with $10$ random restarts.
Since the proposed method performs better with a small value for \textit{sweeps}, we choose the search space $sweeps = 1:3$ and $D = 2:15$.
We summarize the results of the comparison in Table \ref{tab:watertank_comparison}, where the validation runtime includes training and simulation for the whole grid and test runtime indicates the required time of retraining the chosen model structure. 
For the proposed method and the SOTA VTN, the validation runtime includes training and simulation for the whole grid.
\begin{table}[]
    \centering
\begin{tabular}{lcccccc}
\toprule
\multirow{2}{*}{Method} & \multicolumn{2}{c}{Order D} &\multicolumn{2}{c}{Runtime (s)} & RMSE\\
& Grid & Sel. & Val. & Test & Test \\
\midrule
SVDinit & 2,\dots,15 & 11 &\textbf{7} & \textbf{2} & 0.61 \\
VTN \cite{batselier2017tensor} & 2,\dots,7 & 5 &49 & 11 & 0.61 \\
RVS \cite{birpoutsoukis2018efficient} & 1, 2, 3 & 3 & 36180 & 23400 & \textbf{0.54$^*$} \\
LBF \cite{stoddard2018regularized} & 1, 2, 3 & 3& 162 & 120 & 0.56$^*$ \\
\bottomrule
\end{tabular}
\vspace{0.25cm}
    \caption{Comparing runtime and RMSE. SVDinit (proposed method, deterministic initialization), VTN (SOTA VTN) \cite{batselier2017tensor}, RVS \cite{birpoutsoukis2018efficient} and LBF \cite{stoddard2018regularized}. "Sel." denotes selected model structure. For SVDinit, $sweeps = 2$ is chosen. The astrisk $^*$ indicates the use of Bayesian regularization techniques.}
\label{tab:watertank_comparison}
\end{table}\\
\\
The proposed method achieves a $7 \times$ speed-up compared to the SOTA VTN, while maintaining accuracy. 
Interestingly, the choice of the training method has a significant impact on the chosen model order $D$.
Compared to the regularized approaches, both SVDinit and SOTA VTN exhibit a slight RMSE decrease of $8.2\%$ and $11.5\%$ with a more than $23 \times$ and $5100 \times$ speed-up of using SVDinit compared to LBF and RVS, respectively. 
Notable is the achieved $23\times$ increase in efficiency compared to \cite{stoddard2018regularized} despite our method is trained on a weaker hardware, indicating that the algorithmic efficiency gain is even larger than the pure runtime comparison suggests.%

\section{Conclusion}
In this paper, we proposed an efficient automatic algorithm for hyperparameter selection by incrementally increasing the model order $D$ and memory $M$.
We derived it via an LSE formulation for the full Volterra model and transferred it to the VTN format.
In our experiments, we confirmed a significant increase in efficiency with increased or competitive accuracy compared to the literature.
A clear future work direction is to improve the initial guess by e.g. incorporating Bayesian Priors.

\appendix
\section{Appendix}
\subsection{Important Properties and Definitions}\label{app:properties_and_definitions}
\begin{defn}[$k$-Mode Product \cite{kolda2009tensor}]\label{def:kmode_Product}
    The $k$-mode product $\Ct = \At \times_k \B$ of a tensor $\At \in \mathbb{R}^{I_1\times\dots\times I_k \times\dots\times I_D}$ with a matrix $\B \in \mathbb{R}^{J \times I_K}$ is given by
    \begin{align}
    c&\big|_{[i_1,\dots,j,i_{k+1},\dots ,i_D]} = \sum_{i_{k=1}}^{I_k} b\big|_{[j,i_k]} a\big|_{[i_1,\dots,i_k,i_{k+1}\dots,i_D]}, \nonumber
    \end{align}
    such that $\Ct \in \mathbb{R}^{I_1\times\dots\times I_{k-1}\times J \times I_{k+1}\times \dots \times I_D}$.
\end{defn}
The following properties are stated in e.g. \cite{batselier2017tensor, kolda2009tensor}: 
\begin{align}
    &(\A \otimes \B) ^\top = \A^\top \otimes \B^\top, \quad (\A \otimes \B) ^\dagger = \A^\dagger \otimes \B^\dagger \label{eq:kron_prop_0}\\
    &(\A \otimes \B) (\C \otimes \D) = \A\C \otimes \B \D \label{eq:kron_prop_1}\\
    &\vect(\A \B \C^\top) = (\C \otimes \A) \vect{(\B)} \label{eq:kron_prop_2b}\\
    &\Y\big|_{\text{right}} = \U_1 \A\big|_{\text{right}}(\U_3 \otimes \U_2)^\top  \nonumber \dots \\
    &\qquad \qquad \quad\dots \Leftrightarrow \Yt = \At \times_1 \U_1 \times_2\U_2\times_3\U_3 \label{eq:kron_prop_4}
\end{align}
\subsection{Constraint matrix $\C$ and $\C_M$ 
(full VM)} \label{app:constrain_matrix_VM}
We build the constraint matrix $\C \in \mathbb{R}^{MI^D \times I^{D+1}}$ via the binary indicator vector $\cvec^{\otimes D+1} \in \mathbb{R}^{I^{D+1}}$, defined as  
\begin{align}
    c_{[m]}^{\otimes D+1} := \begin{cases}
1,\ & \text{if $\w_{D+1}$ weighs a term created} \\
 & \text{by extending $\U^{\otimes D}$ to $\U^{\otimes D+1}$} \\
0,\ & \text{otherwise}.
    \end{cases} \nonumber 
\end{align}
We construct the rows of $\C$, such that $\cvec_{[m,:]} := \ev_{I^{D+1},m}^\top$.
The matrix $\U^{\otimes D+1}$ is computed based on taking $D+1$ copies of $\uv_n$ with the Kronecker product.
Inserting an additional $\uv_n$ between the $d$th and $d+1$th $\uv_n$ yields
\begin{align}
    \uv_n^{\otimes D+1} &:= \uv_n^{\otimes d} \otimes \uv_n \otimes \uv_n^{\otimes (D-d)}, \label{eqapp:cvec_otimes_perm_1}\\
    \cvec^{\otimes D+1} &:= \mathbf{1}_{I^d} \otimes (\1_{I} - \ev_{I,1}) \otimes \mathbf{1}_{I^{(D-d)}}. \label{eqapp:cvec_otimes_perm_2}
\end{align}
The constraint matrix $\C_M \in \mathbb{R}^{M_{\text{add}}^D\times (I + M_{\text{add}})^D}$ is constructed analogously. 
The main difference is, that we compute the indicator vector $\cvec^{\otimes D}$ as the $D$-times Kronecker product of the vector $\cvec \ := [\mathbf{0}_I \ \mathbf{1}_{M_{\text{add}}}]^\top$.
\subsection{Orthonormal basis $\Z$ and $\Z_M$ (full VM)}\label{app:nullspace_matrix_VM}
We compute $\Z_M \in \mathbb{R}^{(I + M_{\text{add}})^D\times I^D}$ analogue to $\Z$. For $\Z$, we choose the columns of $\Z \in \mathbb{R}^{I^{D+1}\times I^D}$ to span an orthonormal basis for $\mathcal{N}(\C)$.
We construct $\Z$ via the complement $\z^{\otimes D+1} \in \mathbb{R}^{I^{D+1}}$ of the binary indicator vector $\cvec^{\otimes D+1}$.
By replacing the term $(\1_{I} - \ev_{I,1})$ in Equation (\ref{eqapp:cvec_otimes_perm_2}) with its binary complement $\ev_{I,1}$, we obtain 
\begin{align}
    \z^{\otimes D+1} :=&\ \1_{I^d} \otimes \ev_{I,1} \otimes \1_{I^{(D-d)}}. \label{eqapp:zvec_otimes_2}
\end{align}
\subsection{Proof of Lemma \ref{lem:LSU_VM} and Corollary \ref{cor:identical_outputs}} 
\textbf{Forming $\hat{\w}_{\text{LSE}}$ from $\hat{\w}_{\text{LSU}}$:} $\Z$ is a tall matrix constructed with Equation (\ref{eqapp:zvec_otimes_2}). 
For $d=0$, computing $\w_{\text{LSE}} = \Z\w_{\text{LSU}}$ replicates $\w_{\text{LSU}}$ at the position of $\ev_{I,1}$ and fills the remaining entries with zeros as
\begin{align}
    \w_{\text{LSE}} = [\w_{\text{LSU}} \ \mathbf{0}_{I^D}]^\top. \label{eqapp:wLSE_from_wLSU}
\end{align}
\textbf{Equality $\hat{\w}_{\text{LSU}} = \hat{\w}_{D}$:} 
By the construction of $\Z$ as described above, we ensure $\U^{\otimes D+1}\Z = \U^{\otimes D}$.
Inserting this into the LS problem for $\w_{\text{LSU}}$ in Equation (\ref{eq:VM_LSU}) yields a $D$th-order LS problem for $\w_{\text{LSU}}$ equivalent to the $D$th-order LS problem for $\w_{D}$ in Equation (\ref{eq:VM_LS_D}).
Since the minimum-norm solution of a LS problem is unique, we conclude that $\hat{\w}_{\text{LSU}} = \hat{\w}_{D}$. \\
\\
\textbf{The identical input-output behavior} follows from $\hat{\w}_{\text{LSE}} = \Z\hat{\w}_{\text{LSU}} = \Z \hat{\w}_{D}$.
Equation (\ref{eqapp:wLSE_from_wLSU}) shows that the first $I^D$ entries in $\w_{\text{LSE}}$ are equal to $\hat{\w}_{\text{LSU}} = \hat{\w}_{D}$ and weigh the terms that already existed in $\U^{\otimes D}$ and zeroes out the remaining $I^D$ entries in $\U^{\otimes D+1}$.
\subsection{Proof of Lemma \ref{lem:LSU_LSD1_VM}}
We have $\rv_{D+1} \in \mathcal{N}(\U^{\otimes D+1\top})$. Assuming persistently exciting inputs up to order $D+1$ \cite{batselier2017tensor}, $\U^{\otimes D+1}\Z$ is a tall matrix, such that the solution to Equation (\ref{eq:VM_LSU}) is
\begin{align}
    \hat{\w}_{\text{LSU}} &= \underbrace{(\U^{\otimes D+1}\Z)^\dagger}_{=\ (\Z^\top \U^{\otimes D+1\top}\U^{\otimes D+1}\Z)^{-1}\Z^\top\U^{\otimes D+1\top}} \y. \label{eqapp:wLSU_solution}
\end{align}
With $\hat{\w}_{\text{LSE}} = \Z\hat{\w}_{\text{LSU}}$ and Equation (\ref{eqapp:wLSU_solution}), we obtain
\begin{align}
    \hat{\w}_{\text{LSE}} &= \Z \hat{\w}_{\text{LSU}} \nonumber \\
    & = \Z(\Z^\top \U^{\otimes D+1\top}\U^{\otimes D+1}\Z)^{-1}\Z^\top\U^{\otimes D+1\top} \y \nonumber \\
    &=  \Z(\Z^\top \U^*\Z)^{-1}\Z^\top\U^{\otimes D+1\top} \dots \\
    & \quad \dots (\U^{\otimes D+1}\hat{\w}_{D+1} + \rv_{D+1}) \nonumber\\
    &= \underbrace{\Z(\Z^\top \U^{*}\Z)^{-1}\Z^\top\U^*}_{:=\ \Pm} \hat{\w}_{D+1}. \label{eqapp:wLSE_P_wD1}
\end{align} 
$\Pm$ is an oblique ($\U^*$-orthogonal) projector onto $\mathcal{C}(\Z)$ along $\mathcal{N}(\U^{\otimes {D\top}}\U^{\otimes D+1})$ \cite[Theorem 2.2, Theorem 8.1]{stewart2011numerical}.
\subsection{Proof of Corollary \ref{cor:conjugateUpdate_VM}}
Lemma \ref{lem:LSU_LSD1_VM} implies $\hat{\w}_\text{LSE} \in \mathcal{C}(\Z)$. 
With Equations (\ref{eq:updateVector_VM}) and (\ref{eqapp:wLSE_P_wD1}) we express $\hat{\w}_{\Delta} = (\I_{I^{D+1}} - \Pm)\hat{\w}_{D+1}$.
With \cite[Theorem 2.2, Theorem 8.1]{stewart2011numerical}, it follows that $\hat{\w}_{\Delta} \in \mathcal{N}(\U^{\otimes {D\top}}\U^{\otimes D+1})$ and we obtain
\begin{align}
    \hat{\w}_{\text{LSE}}^\top \U^* \hat{\w}_{\Delta} = 0.
\end{align}
\subsection{Proof of Theorem \ref{thm:residualModel_VM} and Theorem \ref{thm:orthogonalUpdate_VM}}
With Corollary \ref{cor:identical_outputs}, Equation (\ref{eq:updateVector_VM}) and (\ref{eq:updateVector_outputs_VM}), we obtain
\begin{align}
    \hat{\y}_{D+1} & = \U^{\otimes D+1}(\hat{\w}_{\text{LSE}} + \hat{\w}_\Delta) = \hat{\y}_{D} + \hat{\y}_\Delta \\
    \hat{\y}_{D}^\top\hat{\y}_{\Delta} &= 0.
\end{align}
To show that $\hat{\w}_\Delta$ is a $D+1$th order VTN model of $\rv_D$, consider $\y = \U^{\otimes D+1}(\hat{\w}_{\text{LSE}} + \hat{\w}_\Delta)  + \rv_{D+1}$ and \\$\y = \U^{\otimes D+1}\hat{\w}_{\text{LSE}} + \rv_D$. It follows \begin{align}
    \U^{\otimes D+1}\hat{\w}_{\text{LSE}} + \rv_D = \U^{\otimes D+1}(\hat{\w}_{\text{LSE}} + \hat{\w}_\Delta)  + \rv_{D+1}. \nonumber 
\end{align}
Solving the above Equation for $\hat{\w}_{\Delta}$ yields
\begin{align}
    \hat{\w}_{\Delta} &= (\U^{\otimes D+1\top}\U^{\otimes D+1})^{-1}\U^{\otimes D+1\top}(\rv_{D} -  \rv_{D+1}) \nonumber \\
    &= (\U^{\otimes D+1\top}\U^{\otimes D+1})^{-1}\U^{\otimes D+1\top} \rv_{D},
\end{align}
since $\rv_{D+1} \in \mathcal{N}(\U^{\otimes D+1\top})$.
It follows immediately, that $\hat{\w}_{\Delta}$ is a residual model, which concludes the proof.  
\subsection{Construction of $\C_{\text{TT}}$ and $\Z_{\text{TT}}$}\label{app:CTTandZTT}
We construct $\C_{TT}$ to select all parameters of $\w^{(d')}_{\text{LSE}}$, that weigh the input terms \mbox{$u(n)\dots u(n-M)$} in $\uv_n$.\\
\underline{For $R_{d'} = 1$: } $\C_{\text{TT}} = \begin{bmatrix} 0 \ \mathbf{1}_{M} \end{bmatrix}$ 
and $\z_{\text{TT}} = \ev_{I,1}$.\\
\underline{For $R_{d'} > 1$}, we obtain the Kronecker product structure
\begin{align}
    \C_{\text{TT}} &= \I_{R_{d'}} \otimes \begin{bmatrix} 0 \ \mathbf{1}_{M} \end{bmatrix} \otimes \I_{R_{d'}}\\
    \Z_{\text{TT}} &= \I_{R_{d'}} \otimes \ev_{I,1} \otimes \I_{R_{d'}}.
\end{align}
\subsection{Construction of $\Qt_Z$}\label{app:constructZTT}
We construct $\Qt_Z^{(d')}$, such that $\Qt_Z^{(d')} \times_2 \uv_n^\top = \I_{R{}_{d'}}$. \underline{For $R_{d'} = 1$:} $ \q_Z^{(d')} = \ev_{I,1}$, analogous to Appendix \ref{app:CTTandZTT}.\\
\underline{For $R_{d'} > 1$} With the $2$-mode product in Definition \ref{def:kmode_Product}, we have to expand $\q_Z^{(d')}$, such that every $i$-th entry turns into its $R_{d'} \times R_{d'}$ matrix equivalent. 
The scalar $1$ turns into $\I_{R_{d'}}$ and each $0$ turns into $\mathbf{0}_{R_{d'}\times R_{d'}}$.
We obtain $\Qt_Z^{(d')}$ in Equation (\ref{eq:incD_Zcore}). With Definition \ref{def:right_left_unfold} 
\begin{align}
    \Q_Z^{(d')}\big|_{\text{left}} &= I_{R_{d'}} \otimes \ev_{I,1} \in \mathbb{R}^{R_{d'}I \times R_{d'}}, \\
    \Q_Z^{(d')}\big|_{\text{right}} &=  \ev_{I,1} \otimes I_{R_{d'}}  \in \mathbb{R}^{R_{d'} \times IR_{d'}}; \label{eqapp:Qz_right}
\end{align}
we get $ \Z_{\text{TT}} = \I_{R_{d'}} \otimes \Q_Z^{(d')} \big|_{\text{left}} =  (\Q_Z^{(d')}\big|_{\text{right}})^\top \otimes I_{R_{d'}}.$ 
\subsection{Proof of Theorem \ref{thm:VTN_obtain_WLSE}}
Assume, we insert $\Qt_{\Z}^{(d')}$ at $d' = D + 1$. With Definition \ref{def:TT} and by construction of $\Qt_{\Z}^{(d')}$, we obtain $\hat{\tilde{\w}}{}^{\text{TT}_{\text{LSE}}}$ by appending zeros to $\hat{\tilde{\w}}{}^{\text{TT}_{\text{LSE}}}$.
As shown in Lemma \ref{lem:LSU_VM}, this is the effect of $\Z$. Shifting $d'$ yields a permutation. 
\subsection{Proof of Theorem \ref{thm:VTN_obtain_WLSE_2} and Corollary \ref{cor:identical_outputs}}\label{app:thm_TN}
First, we consider bringing $\hat{\tilde{\w}}{}^{\text{TT}}_{\text{LSE}}$ into side-$d$-mixed canonical form. 
First, we consider the LQ-decomposition of the right-unfolding $\Wt^{(d)}_D|_{\text{right}}$ of the $d$th TT-core of $\hat{\widetilde{\w}}{}_D^{\text{TT}}$ in Equation (\ref{eq:TT_level_1}).
We define the resulting $\Lm \in \mathbb{R}^{R_{d'} \times R_{d'}}$ and $\Qt^{(d)} \in \mathbb{R}^{R_{d'}\times I \times R_{d'}}$ such that
\begin{align}
   \hat{\widetilde{\w}}{}_D^{(d)} &= \vect(\Qt^{(d)} \times_1 \Lm ) \label{eqapp:Wd_LQ_1}\\ 
    &= \vect(\Lm \Q^{(d')}|_{\text{right}}) \label{eqapp:Wd_LQ_2}\\
    &= \left((\Q^{(d')}|_{\text{right}})^\top \otimes \I_{R_{d'}}\right)\vect (\Lm) \label{eqapp:Wd_LQ_3}, 
\end{align}
We use Definition \ref{def:right_left_unfold} and \ref{def:kmode_Product} to obtain Equation (\ref{eqapp:Wd_LQ_1}) and (\ref{eqapp:Wd_LQ_2}). With Equation (\ref{eq:kron_prop_2b}), we obtain Equation (\ref{eqapp:Wd_LQ_3}). 
Similarly, it follows with Equation (\ref{eqapp:Qz_right}) that 
\begin{align}
     \vect{(\Qt_{Z}^{(d')} \times_3 \Lm})   &= \vect(\Q_Z^{(d')}|_{\text{left}}\ \Lm)  \label{eqapp:Wd_LQz_1}\\
     &=  \left(I_{R_{d'}} \otimes \Q_Z^{(d')}\big|_{\text{left}}\right) \vect{(\Lm)}\\
     &= \Z_{\text{TT}}\vect (\Lm) \label{eqapp:Wd_LQz_3}.
\end{align}
It remains to show that $\hat{\w}_{\text{LSU}}^{(d')} = \Lm$.
With Remark \ref{rem:SISO2MIMO} and since (\ref{eq:kron_prop_0}) also holds true for the row-wise Kronecker product, 
it is sufficient to conduct remainder of the proof for $\hat{y}(n)_D$.
As shown in \cite{batselier2017tensor}, we compute $\hat{y}(n)$
\begin{align}
\hat{y}_D(n) &= \hat{\boldsymbol{\omega}}_{d} \left(\widehat{\Wt}^{(D)} \times_2 \uv_n^\top\right) \hat{\boldsymbol{\omega}}_{d+1} \label{eqapp:eff_pred_1}  \\
&= (\hat{\boldsymbol{\omega}}_{d+1}^\top \otimes \uv_n^\top \otimes \hat{\boldsymbol{\omega}}_{d-1}) \vect(\widehat{\Wt}^{(D)}) \label{eqapp:eff_pred_2}\\
\hat{\boldsymbol{\omega}}_{d-1} &:= (\Qt^{(1)} \times_2 \uv_n^\top) \dots (\Qt^{(d-1)} \times_2 \uv_n^\top) \in \mathbb{R}^{R_{d-1}} \nonumber \\ 
    \hat{\boldsymbol{\omega}}_{d+1} &:= (\Qt^{(d+1)} \times_2 \uv_n^\top) \dots (\Qt^{(D)} \times_2 \uv_n^\top) \in \mathbb{R}^{R_{d}}. \nonumber 
\end{align}
We compare \underline{(a)} solving for $\hat{\w}_{\text{LSU}}^{(d')}$ directly with Equation (\ref{eq:VTN_LSU}); and \underline{(b)} solving for $\hat{\widetilde{\w}}{}_D^{\text{TT}}$ first with $\hat{\boldsymbol{\omega}}_{d-1}$ and $ \hat{\boldsymbol{\omega}}_{d+1}$, followed by using Equation~(\ref{eqapp:Wd_LQ_3}) to compute \\$\vect(\Lm) = \left((\Q_{Z}^{(d')}|_{\text{right}})^\top \otimes \I_{R_{d'}}\right)^\dagger \hat{\widetilde{\w}}{}_D^{(d)}$.\\
\underline{Part (a)}
Analogous to the steps from Equation (\ref{eqapp:Wd_LQz_1}) to Equation (\ref{eqapp:Wd_LQz_3}), we can define $\hat{\W}_{\text{LSU}}^{(d')} \in \mathbb{R}^{R_{d'}\times R_{d'}}$ as $\vect(\hat{\W}_{\text{LSU}}^{(d')}) = \hat{\w}_{\text{LSU}}^{(d')}$ and replace
\begin{align}
    \Z_{\text{TT}}\hat{\w}_{\text{LSU}}^{(d')} = \vect{(\Qt_{Z}^{(d')} \times_3 \hat{\W}_{\text{LSU}}^{(d')}}).
\end{align}
To compute the output $\hat{y}(n)$, we use Equation (\ref{eqapp:eff_pred_1}) and (\ref{eqapp:eff_pred_2}).
Further, with $\boldsymbol{\omega}_{d+1}$ and $\boldsymbol{\omega}_{d-1}$ we define 
\begin{align}
    \boldsymbol{\omega}_{d} := (\Qt^{(d)} \times_2 \uv_n^\top)\boldsymbol{\omega}_{d+1} \label{eqapp:omega_d}
\end{align} 
and by construction of $\Qt_{Z}^{(d')}$ obtain
\begin{align}
    \hat{y}(n) &= \boldsymbol{\omega}_{d-1} \left(\Qt_{Z}^{(d')} \times_3 \hat{\W}_{\text{LSU}}^{(d')}\times_2 \uv_n^\top\right) \boldsymbol{\omega}_{d},\\
    &= \boldsymbol{\omega}_{d-1} \left(\I_{R_{d'}}\ \hat{\W}_{\text{LSU}}^{(d')}\right) \boldsymbol{\omega}_{d} \\
    & = \left(\boldsymbol{\omega}_{d}^\top \otimes \boldsymbol{\omega}_{d-1} \right) \vect( \hat{\W}_{\text{LSU}}^{(d')}) \label{eqapp:y_LSU}.
\end{align}
\underline{Part (b)}
Inserting Equation (\ref{eqapp:Wd_LQ_3}) in Equation (\ref{eqapp:eff_pred_2}) yields
\begin{align} 
&\hat{y}(n) = (\hat{\boldsymbol{\omega}}_{d+1}^\top \otimes \uv_n^\top \otimes \hat{\boldsymbol{\omega}}_{d-1})\dots \nonumber \\
&\quad \dots \left((\Q^{(d)}|_{\text{right}})^\top \otimes \I_{R_{d'}}\right)\vect (\Lm) \nonumber\\
&= \left( \left[(\hat{\boldsymbol{\omega}}_{d+1}^\top \otimes \uv_n^\top)(\Q^{(d)}|_{\text{right}})^\top\right] \otimes \hat{\boldsymbol{\omega}}_{d-1} \right) \vect (\Lm) \nonumber\\
&= \left( \left[\left(\Q^{(d)}|_{\text{right}})(\hat{\boldsymbol{\omega}}_{d+1}^\top \otimes \uv_n^\top\right)^\top\right]^\top \otimes \hat{\boldsymbol{\omega}}_{d-1} \right) \vect (\Lm) \nonumber\\
&= \left( \left[(\Qt^{(d)} \times_2 \uv_n^\top) \hat{\boldsymbol{\omega}}_{d+1}\right]^\top \otimes \hat{\boldsymbol{\omega}}_{d-1} \right) \vect (\Lm) \nonumber\\
&= \left( \boldsymbol{\omega}_d^\top \otimes \hat{\boldsymbol{\omega}}_{d-1} \right) \vect (\Lm), \label{eqapp:y_L}
\end{align}
where we obtain the second line with Equation (\ref{eq:kron_prop_1}), the third line with Equation (\ref{eq:kron_prop_0}), the fourth line with Equation (\ref{eq:kron_prop_4}) and Definition \ref{def:kmode_Product}, and the fifth and last line with Equation (\ref{eqapp:omega_d}).
With the equality of Equation (\ref{eqapp:y_LSU}) and Equation (\ref{eqapp:y_L}), it follows that solving for $\vect( \hat{\W}_{\text{LSU}}^{(d')})$ and for $\vect(\Lm)$ yields identical results, thus $\vect(\hat{\W}_{\text{LSU}}^{(d')}) = \vect(\Lm)$, which concludes the proof.

\bibliographystyle{plain}        
\bibliography{main.bib}           

\end{document}